\documentclass{article}%
\usepackage{amsmath}
\usepackage{graphicx}%
\usepackage{amsfonts}%
\usepackage{amssymb}
\newtheorem{theorem}{Theorem}

\newtheorem{algorithm}[theorem]{Algorithm}

\newtheorem{corollary}[theorem]{Corollary}

\newtheorem{definition}[theorem]{Definition}
\newtheorem{example}[theorem]{Example}

\newtheorem{lemma}[theorem]{Lemma}

\newtheorem{proposition}[theorem]{Proposition}
\newtheorem{remark}[theorem]{Remark}

\newenvironment{proof}[1][Proof]{\textbf{#1.} }{\ \rule{0.5em}{0.5em}}

\begin{document}

\title{Foliating Metric Spaces:\\{\normalsize A Generalization of Frobenius' Theorem}}
\author{Craig Calcaterra}
\maketitle

\begin{abstract}
Using families of curves to generalize vector fields, the Lie bracket is
defined on a metric space, $M$. For $M$ complete, versions of the local and
global Frobenius theorems hold, and flows are shown to commute if and only if
their bracket is zero.

An example is given showing $L^{2}\left(  \mathbb{R}\right)  $ is controllable
by two elementary flows.



\end{abstract}

\section{Introduction}

The main goal of this monograph is to further the point of view that many
beautiful geometrical and analytical results valid on differentiable manifolds
hold on general metric spaces. Besides the wider relevance gained by
generalization, the foundations of the subject are clarified when the limits
of applicability are explored. This effort has a long and often disjointed
history, only one sliver of which is relevant here. The approach in this
paper, which has been used by several others, is to use the well-known
characterization of a vector in a tangent space as an equivalence class of
curves which are tangent to each other. A \textbf{curve} $c$ on a metric space
$\left(  M,d\right)  $ is a continuous map $c:\left(  \alpha,\beta\right)
\rightarrow M$ where $\left(  \alpha,\beta\right)  \subset\mathbb{R}$. Two
curves $c_{i}:\left(  \alpha_{i},\beta_{i}\right)  \rightarrow M$ for $i=1,2$
are \textbf{tangent at }$t\in\left(  \alpha_{1},\beta_{1}\right)  \cap\left(
\alpha_{2},\beta_{2}\right)  $ if%
\[
\underset{h\rightarrow0}{\lim}\frac{d\left(  c_{1}\left(  t+h\right)
,c_{2}\left(  t+h\right)  \right)  }{h}=0\text{.}%
\]
In this way we may generalize a vector field (a family of vectors) on a
manifold as an arc field (a family of curves) on a metric space--Definition
\ref{ArcField}, below.

It has been said the three pillars of differential geometry are: $\left(
I\right)  $ the Inverse Function Theorem, $\left(  II\right)  $ the Existence
Theorem for ordinary differential equations (ODEs) and $\left(  III\right)  $
Frobenius' Theorem. All of these classical theorems may be written with vector
fields on manifolds and so may also be written with arc fields on metric
spaces. We expect any result on manifolds which has a sufficiently
geometrically realized proof can be generalized to metric spaces using curves
in place of vectors. A metric space version of $\left(  I\right)  $ is
contained in \cite{Aubin2}, e.g.; and versions of $\left(  II\right)  $ have
been proven several times independently in e.g., \cite{Panasyuk2},
\cite{Aubin2}, and \cite{CalcBleecker}--see Theorem \ref{CL} below. A version
of $\left(  III\right)  $ is the main result of this paper, Theorem
\ref{FrobeniusThm}: an involutive distribution on a complete metric space is
integrable. Since the result is for complete metric spaces, it generalizes the
classical result on Banach manifolds (proven, e.g., in \cite{abrahamMarsden}).
Theorem \ref{FrobeniusThm} further generalizes the classical result by
assuming only Lipschitz-type regularity instead of smoothness, which is of
interest in, for example, control theory.

As far as I have been able to determine, this particular approach to the proof
of Frobenius' classical theorem has not been vetted in the literature--though
it uses basic, well-known ideas. We outline the approach in this paragraph,
simplified to vector fields on a manifold. The terminology and assumptions
will be clarified in the main body of the paper, and Figures
\ref{FigFrobProof1} and \ref{FigFrobProof2} from Section \ref{SectionFrobThm}
may aid intuition. The crux of the local Frobenius result in two dimensions is
as follows: Given two transverse vector fields $f,g:M\rightarrow TM$ there
exists an integral surface (tangent to linear combinations of $f$ and $g$)
through any point $x_{0}\in M$ when the Lie bracket satisfies $\left[
f,g\right]  =af+bg$ for some choice of functions $a,b:M\rightarrow\mathbb{R}$
(involutivity of $f$ and $g$). To prove this, define
\[
S:=\left\{  F_{t}G_{s}\left(  x_{0}\right)  \in M:\left|  s\right|  ,\left|
t\right|  <\delta\right\}
\]
where $F$ and $G$ are the local flows of $f$ and $g$ guaranteed to exist by
$\left(  II\right)  $. Since $f$ and $g$ are transverse, we may choose
$\delta>0$ small enough for $S$ to be a well-defined surface. $S$ will be
shown to be the desired integral surface through $x_{0}$. Notice $S$ is
tangent to $f$ by construction, but it is not immediately clear $S$ is tangent
to $a^{\prime}f+b^{\prime}g$ for arbitrarily chosen $a^{\prime},b^{\prime}%
\in\mathbb{R}$. Notice, though, that by construction $S$ is tangent to $g$ at
any point $x=G_{s}\left(  x_{0}\right)  $, and also $S$ is tangent to
$a^{\prime\prime}f+b^{\prime\prime}g$ at $x$ for functions $a^{\prime\prime}$
and $b^{\prime\prime}$. Therefore establishing%
\begin{equation}
\left(  F_{t}\right)  ^{\ast}\left(  a^{\prime}f+b^{\prime}g\right)
=a^{\prime\prime}f+b^{\prime\prime}g\text{\quad at\quad}x=G_{s}\left(
x_{0}\right)  \label{FrobOutline5}%
\end{equation}
for some functions $a^{\prime\prime}$ and $b^{\prime\prime}$, proves $S$ is
tangent to $a^{\prime}F+b^{\prime}G$ at an arbitrary point $z=F_{t}%
G_{s}\left(  x_{0}\right)  \in S$, since the push-forward $\left(
F_{t}\right)  _{\ast}$ and the pull-back $\left(  F_{t}\right)  ^{\ast}$ are
inverse to each other and preserve tangency since they are local
lipeomorphisms. Next since the Lie bracket equals the Lie derivative,%
\[
\underset{h\rightarrow0}{\lim}\frac{F_{h}^{\ast}\left(  g\right)  -g}%
{h}=\left[  f,g\right]  =af+bg
\]
by involutivity so%
\[
F_{h}^{\ast}\left(  g\right)  =g+h\left(  af+bg\right)  +o\left(  h\right)
=\widetilde{a}f+\widetilde{b}g+o\left(  h\right)  \text{.}%
\]
Using the fact that $F_{h}^{\ast}\left(  f\right)  =f$ for any $h$, and the
linearity of pullback for fixed $t$, we have for functions $a_{i}$ and
$b_{i}:M\rightarrow\mathbb{R}$%
\[
F_{t/n}^{\ast}\left(  a_{i}f+b_{i}g\right)  =\left(  a_{i+1}f+b_{i+1}g\right)
+o\left(  1/n\right)
\]
for some functions $a_{i+1}$ and $b_{i+1}$. Then since%
\[
F_{t}^{\ast}=\underset{n\text{ times}}{\underbrace{F_{t/n}^{\ast}F_{t/n}%
^{\ast}...F_{t/n}^{\ast}}}=\left(  F_{t/n}^{\ast}\right)  ^{\left(  n\right)
}%
\]
(where the superscript in round brackets denotes composition $n$ times) we get
$\left(  \text{\ref{FrobOutline5}}\right)  $ as follows:%
\begin{align*}
F_{t}^{\ast}\left(  a_{0}f+b_{0}g\right)   &  =\underset{n\rightarrow\infty
}{\lim}\left(  F_{t/n}^{\ast}\right)  ^{\left(  n\right)  }\left(
a_{0}f+b_{0}g\right)  \\
&  =\underset{n\rightarrow\infty}{\lim}a_{n}f+b_{n}g+no\left(  1/n\right)
=a_{\infty}f+b_{\infty}g+0
\end{align*}
completing the sketch for manifolds.

A pivotal fact on which the metric space version relies is that arc fields
which satisfy certain Lipschitz-type conditions generate unique local flows
(proven in \cite{CalcBleecker} and reviewed in Section \ref{SectionReview}).
Also a natural linear structure may be associated with a metric space (though
it has no \textit{a priori} linear structure) using compositions of flows
which faithfully generalizes the linearity of vector fields; this was
introduced in \cite{Colombo}. We present this in Section
\ref{SectionBracket&linearity} along with the generalization of the Lie
bracket for vector fields which uses the well-known asymptotic
characterization of the Lie bracket; i.e., successively follow the flows
forward and backward for time $\sqrt{t}$. This investigation further clarifies
for us the surprising fact Sussman and others have noted: smoothness is not
necessary to define a geometrically meaningful Lie bracket. In Section
\ref{PushFwdSect}, the pull-back along a flow is shown to behave naturally
with linearity and the bracket, which mimics properties of the Lie derivative
on manifolds. Many more such algebraic properties are valid than are contained
in these sections, but in this monograph we present only the minimum machinery
directly relevant to proving Frobenius' Theorem in Section
\ref{SectionFrobThm}.

Section \ref{SectionDist&Foliations} applies this local Frobenius theorem to
study foliations yielding a global theorem on metric spaces. A metric space
generalization of the Nagumo-Brezis Invariance Theorem is proven, which is
used to show integrable distributions are involutive. We do not discuss the
facet of the classical Global Frobenius Theorem which guarantees local
coordinates on which there exist coordinate vector fields tangent or
perpendicular to an involutive distribution. In light of these results,
however, this now seems ripe for exploration.

Section \ref{SectionCommut} proves a well-known result from Hamiltonian
dynamics is also valid for metric spaces: two flows commute if and only if the
bracket is $0$. This is not exactly a corollary of the metric space Frobenius
Theorem, but the proof is a mere simplification of that from Theorem
\ref{FrobeniusThm}.

Finally in Section \ref{SectionExs} an almost trivial example applying these
ideas has a result which astounded me: Any Lebesgue square-integrable function
may be approximated using successive compositions of two elementary flows,
starting from the constant zero function. In other words, $L^{2}\left(
\mathbb{R}\right)  $ is controllable by two flows. You may skip straight to
this Example \ref{ExL2decomp} after perusing the following review and the
definitions in Section \ref{SectionBracket&linearity}. \cite{Sontag} is an
accessible text introducing the terminology of control theory with remarks and
references on infinite dimensional controllability.

\section{\label{SectionReview}Review of terminology and basic results}

The proofs of all of the results from this section are contained in
\cite{CalcBleecker} for forward flows, also called semi-flows. Minimal
changes, stated here, give us the corresponding results for (bidirectional) flows.

A \textbf{metric space} $\left(  M,d\right)  $ is a set of points $M$ with a
function $d:M\times M\rightarrow\mathbb{R}$ called the \textbf{metric} which
has the following properties:%
\[%
\begin{array}
[c]{lll}%
\left(  \text{i}\right)  & d(x,y)\geq0 & \text{\textbf{positivity}}\\
\left(  \text{ii}\right)  & d(x,y)=0\text{\quad\textit{iff\quad}}x=y &
\text{\textbf{nondegeneracy}}\\
\left(  \text{iii}\right)  & d(x,y)=d(y,x) & \text{\textbf{symmetry}}\\
\left(  \text{iv}\right)  & d(x,y)\leq d(x,z)+d(z,y) & \text{\textbf{triangle
inequality}}%
\end{array}
\]
for all $x,y,z\in M$. The open ball of radius $r$ about $x\in M$ is denoted by
$B\left(  x,r\right)  :=\left\{  y:d\left(  x,y\right)  <r\right\}  $. We
assume throughout this paper that $\left(  M,d\right)  $ is a locally complete
metric space, i.e., there exists a complete neighborhood of each point in $M$.
Denote the open ball in $M$ about $x_{0}\in M$ with radius $r$ by%
\[
B\left(  x_{0},r\right)  :=\left\{  x\in M:d\left(  x,x_{0}\right)
<r\right\}  \text{.}%
\]
A map $f:\left(  M,d_{M}\right)  \rightarrow\left(  N,d_{N}\right)  $ between
metric spaces is \textbf{Lipschitz} continuous if there exists $K_{f}\geq0$
such that%
\[
d_{N}\left(  f\left(  x_{1}\right)  ,f\left(  x_{2}\right)  \right)  \leq
K_{f}d_{M}\left(  x_{1},x_{2}\right)
\]
for all $x_{1},x_{2}\in M$. A \textbf{lipeomorphism} is an invertible
Lipschitz map whose inverse is also Lipschitz (i.e., stronger than a
homeomorphism, weaker than a diffeomorphism).

The following definition is made in analogy with vector fields on manifolds,
where vectors are represented as curves on the manifold.

\begin{definition}
\label{ArcField}An \textbf{arc field} on $M$ is a continuous map
$X:M\times\left[  -1,1\right]  \rightarrow M$ such that for all $x\in M$,
$X\left(  x,0\right)  =x$,%
\[
\rho\left(  x\right)  :=\sup_{s\neq t}\frac{d\left(  X\left(  x,s\right)
,X\left(  x,t\right)  \right)  }{\left|  s-t\right|  }<\infty,
\]
$($i.e., $X\left(  x,\cdot\right)  $ is Lipschitz$)$, and the function
$\rho\left(  x\right)  $ is locally bounded so%
\[
\rho\left(  x,r\right)  :=\sup_{y\in B\left(  x,r\right)  }\left\{
\rho\left(  y\right)  \right\}  <\infty,
\]
for $r>0$ sufficiently small.

A \textbf{solution curve }to $X$ is a curve $\sigma$ tangent to $X$, i.e.,
$\sigma:\left(  \alpha,\beta\right)  \rightarrow M$ for some open interval
$\left(  \alpha,\beta\right)  \subset\mathbb{R}$ has the following property
for each $t\in\left(  \alpha,\beta\right)  $%
\begin{equation}
\lim_{h\rightarrow0}\frac{d\left(  \sigma\left(  t+h\right)  ,X\left(
\sigma\left(  t\right)  ,h\right)  \right)  }{h}=0\text{,} \label{solncond1}%
\end{equation}
i.e., $d\left(  \sigma\left(  t+h\right)  ,X\left(  \sigma\left(  t\right)
,h\right)  \right)  =o\left(  h\right)  $.
\end{definition}

$\rho$ is a bound on the speed of the arcs. $\alpha$ and $\beta$ are members
of the extended reals $\mathbb{R}\cup\left\{  \pm\infty\right\}  $.

The two variables for arc fields and flows which are usually denoted by $x$
and $t$ are often thought of as representing space and time. In this paper
$x,y,$ and $z$ are used for space variables, while $r,s,t,$ and $h$ may fill
the time variable slot. An arc field $X$ will often have its variables migrate
liberally between parentheses and subscripts%
\[
X\left(  x,t\right)  =X_{x}\left(  t\right)  =X_{t}\left(  x\right)
\]
depending on which variable we wish to emphasize in a calculation. We also use
this convention for flows $F$ defined below.

The following conditions guarantee existence and uniqueness of solutions.

\bigskip\noindent\textbf{Condition E1:} For each $x_{0}\in M$, there are
positive constants $r,\delta$ and $\Lambda_{X}$ such that for all $x,y\in
B\left(  x_{0},r\right)  $ and $t\in\left(  -\delta,\delta\right)  $%
\[
d\left(  X_{t}\left(  x\right)  ,X_{t}\left(  y\right)  \right)  \leq d\left(
x,y\right)  \left(  1+\left|  t\right|  \Lambda_{X}\right)  \text{.}%
\]

\noindent\textbf{Condition E2}:%
\[
d\left(  X_{s+t}\left(  x\right)  ,X_{t}\left(  X_{s}\left(  x\right)
\right)  \right)  =O\left(  st\right)
\]
as $st\rightarrow0$ locally uniformly in $x$; in other words, for each
$x_{0}\in M$, there are positive constants $r,\delta$ and $\Omega_{X}$ such
that for all $x\in B\left(  x_{0},r\right)  $ and $s,t\in\left(
-\delta,\delta\right)  $%
\[
d\left(  X_{s+t}\left(  x\right)  ,X_{t}\left(  X_{s}\left(  x\right)
\right)  \right)  \leq\left|  st\right|  \Omega_{X}\text{.}%
\]
\bigskip%

%TCIMACRO{\FRAME{ftbphFU}{4.1381in}{1.8161in}{0pt}{\Qcb{Conditions E1 and E2}%
%}{}{conde1ande2.eps}{\special{ language "Scientific Word";  type "GRAPHIC";
%width 4.1381in;  height 1.8161in;  depth 0pt;  original-width 3.5993in;
%original-height 1.5645in;  cropleft "0";  croptop "1";  cropright "1";
%cropbottom "0";  filename 'CondE1andE2.eps';file-properties "XNPEU";}} }%
%BeginExpansion
\begin{figure}
[ptbh]
\begin{center}
\includegraphics[
height=1.8161in,
width=4.1381in
]%
{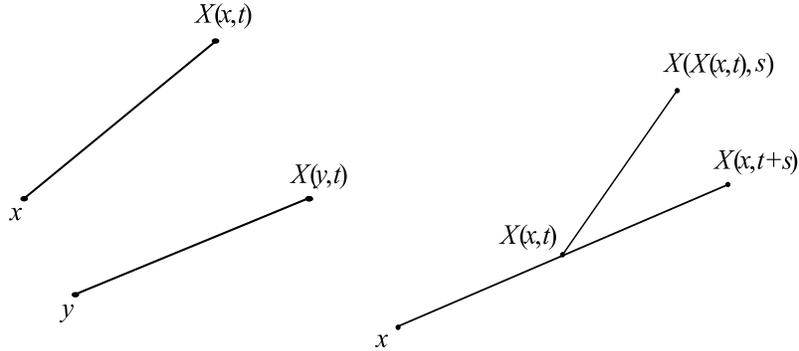}%
\caption{Conditions E1 and E2}%
\end{center}
\end{figure}
%EndExpansion

\begin{theorem}
\label{CL}Let $X$ be an arc field satisfying E1 and E2 on a locally complete
metric space $M$. Then given any point $x\in M,$ there exists a unique
solution $\sigma_{x}:\left(  \alpha_{x},\beta_{x}\right)  \rightarrow M$ with
$\sigma_{x}\left(  0\right)  =x$.
\end{theorem}

Several remarks are in order. Here, $x$ is called the \textbf{initial
condition} for the solution $\sigma_{x}$ in the above theorem. Uniqueness of
solutions means that for any $x\in M$, the curve $\sigma_{x}:\left(
\alpha_{x},\beta_{x}\right)  \rightarrow M$ has maximal domain $\left(
\alpha_{x},\beta_{x}\right)  $ in the sense that for any other solution
$\widehat{\sigma}_{x}:\left(  \widehat{\alpha}_{x},\widehat{\beta}_{x}\right)
\rightarrow M$ also having initial condition $x$, we have $\left(
\widehat{\alpha}_{x},\widehat{\beta}_{x}\right)  \subset\left(  \alpha
_{x},\beta_{x}\right)  $ and $\widehat{\sigma}_{x}=\sigma_{x}|_{\left(
\widehat{\alpha}_{x},\widehat{\beta}_{x}\right)  }$ $($i.e., $\sigma_{x}$ is
the \textbf{maximal solution curve}$)$.

The proof of Theorem \ref{CL} is constructive and shows the \textbf{Euler
curves} $X_{t/n}^{\left(  n\right)  }\left(  x\right)  $ converge to the
solution. Here we are using $f^{\left(  n\right)  }$ to denote the composition
of a map $f:M\rightarrow M$ with itself $n$ times so%
\[
X_{\frac{t}{n}}^{\left(  n\right)  }\left(  x\right)  =\text{ }\underset
{n\text{ times}}{\underbrace{X_{\frac{t}{n}}\circ X_{\frac{t}{n}}\circ...\circ
X_{\frac{t}{n}}}}\left(  x\right)
\]
and we have
\[
\underset{n\rightarrow\infty}{\lim}X_{\frac{t}{n}}^{\left(  n\right)  }\left(
x\right)  =\sigma_{x}\left(  t\right)  \text{.}%
\]
for suitably small $\left|  t\right|  $. Other, slightly different
formulations of Euler curves also lead to the same result, $\sigma$, under
Conditions E1 and E2, e.g.,%
\[
\xi_{n}\left(  t\right)  :=X_{t-i\cdot2^{-n}}X_{2^{-n}}^{\left(  i\right)
}\left(  x\right)  \text{\quad for \quad}i\cdot2^{-n}\leq t\leq\left(
i+1\right)  2^{-n}%
\]
also has%
\[
\underset{n\rightarrow\infty}{\lim}\xi_{n}\left(  t\right)  =\sigma_{x}\left(
t\right)
\]
for suitably small $\left|  t\right|  $.

Theorem \ref{CL} and those that follow are true under more general conditions
outlined in \cite{CalcBleecker} and \cite{Panasyuk2}. But throughout this
paper and in every application I've seen, E1 and E2 are satisfied and are
\underline{E}asier to use.

\begin{example}
\label{Banach Example}A \textbf{Banach space} $\left(  M,\left\|
\cdot\right\|  \right)  $ is a normed vector space, complete in its norm
$($e.g., $\mathbb{R}^{n}$ with Euclidean norm$)$. A Banach space is an example
of a metric space with $d\left(  u,v\right)  :=\left\|  u-v\right\|  $. A
\textbf{vector field} on a Banach space $M$ is a map $f:M\rightarrow M$. A
\textbf{solution} to a vector field $f$ with \textbf{initial condition} $x$ is
a curve $\sigma_{x}:\left(  \alpha,\beta\right)  \rightarrow M$ defined on an
open interval $\left(  \alpha,\beta\right)  \subset\mathbb{R}$ containing $0 $
such that $\sigma_{x}\left(  0\right)  =x$ and $\sigma_{x}^{\prime}\left(
t\right)  =f\left(  \sigma_{x}\left(  t\right)  \right)  $ for all
$t\in\left(  \alpha,\beta\right)  $. The classical Picard-Lindel\"{o}f Theorem
guarantees unique solutions for any locally Lipschitz $f$. With a few tricks,
most differential equations can be represented as vector fields on a suitably
abstract space.

Every Lipschitz vector field $f:M\rightarrow M$ gives rise to an arc field
$X\left(  x,t\right)  :=x+tf\left(  x\right)  $ and it is easy to check $X$
satisfies E1 and E2 $($cf. \cite{CalcBleecker}$)$. Further the solutions to
the arc field are exactly the solutions to the vector field. Therefore Theorem
\ref{CL} is a generalization of the classical Picard-Lindel\"{o}f Theorem.
\end{example}

\begin{remark}
\label{RemUniformSolutions}Of prime import for this monograph, the proof of
Theorem \ref{CL} actually shows solutions are \textbf{locally uniformly
tangent} to $X$:
\[
d\left(  X_{x}\left(  t\right)  ,\sigma_{x}\left(  t\right)  \right)
=o\left(  t\right)
\]
locally uniformly for $x\in M$, i.e., for each $x_{0}\in M$ there exists a
constant $r>0$ such that for any $\varepsilon>0$ there exists a $\delta>0$
such that for all $x\in B\left(  x_{0},r\right)  $%
\[
\frac{d\left(  X_{x}\left(  t\right)  ,\sigma_{x}\left(  t\right)  \right)
}{\left|  t\right|  }<\varepsilon
\]
whenever $0<\left|  t\right|  <\delta$.

More than this, the proof also shows solutions are tangent uniformly for all
arc fields $X$ which satisfy $E1$ and $E2$ for specified $\Lambda$ and
$\Omega$.
\end{remark}

We denote local uniform tangency of two arc fields $X$ and $Y$ by $X\sim Y$.
It is easy to check $\sim$ is an equivalence relation. E.g., transitivity
follows from the triangle inequality:%
\[
\frac{d\left(  X_{t}\left(  x\right)  ,Z_{t}\left(  x\right)  \right)
}{\left|  t\right|  }\leq\frac{d\left(  X_{t}\left(  x\right)  ,Y_{t}\left(
x\right)  \right)  }{\left|  t\right|  }+\frac{d\left(  Y_{t}\left(  x\right)
,Z_{t}\left(  x\right)  \right)  }{\left|  t\right|  }\text{.}%
\]
We use the symbol $\sim$ in many contexts in this paper (particularly Section
\ref{SectionDist&Foliations}), but there is always a local-uniform-tangency
property associated with it.

\begin{corollary}
Assume the conditions of Theorem \ref{CL} and let $s\in\left(  \alpha
_{x},\beta_{x}\right)  $ and $y=\sigma_{x}(s)$. \ Then $\alpha_{y}=\alpha
_{x}-s$ and $\beta_{y}=\beta_{x}-s$ so%
\[
\left(  \alpha_{y},\beta_{y}\right)  =\left(  \alpha_{\sigma_{x}(s)}%
,\beta_{\sigma_{x}(s)}\right)  =\{t:\alpha_{x}-s<t<\beta_{x}-s\}\text{.}%
\]
Thus $t\in\left(  \alpha_{y},\beta_{y}\right)  $ if and only if $t+s\in\left(
\alpha_{x},\beta_{x}\right)  ,$ and then we have
\[
\sigma_{\sigma_{x}(s)}(t)=\sigma_{x}(s+t)\text{.}%
\]
Defining $W\subset M\times\mathbb{R}$ by%
\[
W:=\{(x,t)\in M\times\mathbb{R}:t\in\left(  \alpha_{x},\beta_{x}\right)  \}
\]
and $F:W\rightarrow M$ by $F(x,t):=\sigma_{x}(t)$ we have:\bigskip

\noindent$\left(  i\right)  $\quad$M\times\{0\}\subset W$ and $F(x,0)=x$ for
all $x\in M$.

\noindent$\left(  ii\right)  $\quad For each (fixed) $x\in M$, $F(x,\cdot
):\left(  \alpha_{x},\beta_{x}\right)  \rightarrow M$ is the maximal solution
$\sigma_{x}$ to $X$.

\noindent$\left(  iii\right)  $\quad$F(t,F(s,x))=F(t+s,x)$.\bigskip
\end{corollary}

$F$ is called the \textbf{local flow} generated by the arc field $X$. Compare
Condition E2 with property $\left(  iii\right)  $ above to see why an arc
field might be thought of as a ``pre-flow''.

\begin{theorem}
\label{ExpGrowth}Let $\sigma_{x}:\left(  \alpha_{x},\beta_{x}\right)
\rightarrow M$ and $\sigma_{y}:\left(  \alpha_{y},\beta_{y}\right)
\rightarrow M$ be two solutions to an arc field $X$ which satisfies E1. Assume
$\left(  \alpha_{x},\beta_{x}\right)  \cap\left(  \alpha_{y},\beta_{y}\right)
\supset I$ for some interval $I$, and assume $\Lambda_{X}$ holds on a set
containing%
\[
\left\{  \sigma_{x}\left(  t\right)  :t\in I\right\}  \cup\left\{  \sigma
_{y}\left(  t\right)  :t\in I\right\}  \text{.}%
\]
Then
\[
d\left(  \sigma_{x}\left(  t\right)  ,\sigma_{y}\left(  t\right)  \right)
\leq e^{\Lambda_{X}\left|  t\right|  }d\left(  x,y\right)  \text{\ for all
}t\in I\text{,}%
\]
i.e.,%
\begin{equation}
d\left(  F_{t}\left(  x\right)  ,F_{t}\left(  y\right)  \right)  \leq
e^{\Lambda_{X}\left|  t\right|  }d\left(  x,y\right)  \text{.}
\label{CondExpGrowth}%
\end{equation}
\end{theorem}

\begin{theorem}
For $F$ and $W$ as above, $W$ is open in $M\times\mathbb{R}$ and $F$ is
continuous on $W$.
\end{theorem}

For fixed $t$ it is clear $F_{t}$ is a local lipeomorphism, when defined, by
Theorem \ref{ExpGrowth}. Compare Condition E1 with line $\left(
\text{\ref{CondExpGrowth}}\right)  $ to see why E1 may be thought of as a
local linearity property for $X$, needed for the continuity of $F$.

\begin{definition}
\label{LinSpeed}An arc field $X$ on a metric space $M$ is said to have
\textbf{linear speed growth} if there is a point $x\in M$ and positive
constants $c_{1}$ and $c_{2}$ such that for all $r>0$
\begin{equation}
\rho\left(  x,r\right)  \leq c_{1}r+c_{2}, \label{LinGro}%
\end{equation}
where $\rho\left(  x,r\right)  $ is the local bound on speed given in
Definition \ref{ArcField}.
\end{definition}

\begin{theorem}
\label{LongTime}Let $X$ be an arc field on a complete metric space $M$, which
satisfies E1 and E2 and has linear speed growth. Then $F$ is a $($full$)$
\textbf{flow} with domain $W=M\times\mathbb{R}$.
\end{theorem}

\begin{example}
\label{Flow<-->ArcField}Every local flow on a metric space is generated by an
arc field. Any local flow $F$ gives rise to an arc field $\overline{F}%
:M\times\left[  -1,1\right]  \rightarrow M$ defined by%
\[
\overline{F}\left(  x,t\right)  :=\left\{
\begin{array}
[c]{ll}%
F\left(  x,t\right)  & \text{if }t\in\left(  \frac{\alpha_{x}}{2}%
,\frac{\beta_{x}}{2}\right) \\
F\left(  x,\frac{\alpha_{x}}{2}\right)  & \text{if }t\in\left[
-1,\frac{\alpha_{x}}{2}\right] \\
F\left(  x,\frac{\beta_{x}}{2}\right)  & \text{if }t\in\left[  \frac{\beta
_{x}}{2},1\right]  \text{.}%
\end{array}
\right.
\]
$($The issue here is that $F$, being a \textit{local} flow, may have
$\alpha_{x}$ or $\beta_{x}<1$.$)$ Clearly the local flow generated by
$\overline{F} $ is $F$. Since all our concerns with arc fields are local, we
will never focus on $t\notin\left(  \frac{\alpha_{x}}{2},\frac{\beta_{x}}%
{2}\right)  $ and henceforth we will not notationally distinguish between
$\overline{F}$ and $F$ as arc fields.
\end{example}

With this identification of flows being arc fields (but not necessarily
\textit{vice}-\textit{versa}) we may simplify Remark \ref{RemUniformSolutions}
to: $X\sim F$ if $X$ satisfies E1 and E2.

\section{The bracket and linearity\label{SectionBracket&linearity}}

To simplify notation we drop parentheses for expressions such as $Y_{t}\circ
X_{s}\left(  x\right)  =Y_{t}\left(  X_{s}\left(  x\right)  \right)  $ and
write $Y_{t}X_{s}\left(  x\right)  $ since the composition of arbitrary maps
is associative.

\begin{definition}
The \textbf{bracket} of two arc fields $X$ and $Y$ is the map $\left[
X,Y\right]  :M\times\left[  -1,1\right]  \rightarrow M$ with%
\begin{equation}
\left[  X,Y\right]  \left(  x,t\right)  :=\left\{
\begin{array}
[c]{c}%
Y_{-\sqrt{t}}X_{-\sqrt{t}}Y_{\sqrt{t}}X_{\sqrt{t}}\left(  x\right) \\
X_{-\sqrt{\left|  t\right|  }}Y_{-\sqrt{\left|  t\right|  }}X_{\sqrt{\left|
t\right|  }}Y_{\sqrt{\left|  t\right|  }}\left(  x\right)
\end{array}
\right.
\begin{array}
[c]{c}%
\text{for }t\geq0\\
\text{for }t<0\text{.}%
\end{array}
\label{bracketDef}%
\end{equation}
\end{definition}

There are many different equivalent characterizations of the Lie bracket on a
manifold. $\left(  \ref{bracketDef}\right)  $ uses the obvious choice of the
\textbf{asymptotic} characterization to generalize the concept to metric
spaces. $\left[  X,Y\right]  \left(  x,t\right)  $ traces out a small
``parallelogram'' in $M$ starting at $x$, which hopefully almost returns to
$x$. The bracket measures the failure of $X$ and $Y$ to commute as will be
made clear in Theorems \ref{ThmCommute} and \ref{FrobeniusThm}.

\begin{definition}
We say $X$ \textbf{\&} $Y$ \textbf{close} if%
\[
d\left(  Y_{s}X_{t}\left(  x\right)  ,X_{t}Y_{s}\left(  x\right)  \right)
=O\left(  \left|  st\right|  \right)
\]
locally uniformly in $x$, i.e., if for each $x_{0}\in M$ there exist positive
constants $C_{XY},\delta,$ and $r$ such that for all $x\in B\left(
x_{0},r\right)  $%
\[
d\left(  Y_{s}X_{t}\left(  x\right)  ,X_{t}Y_{s}\left(  x\right)  \right)
\leq C_{XY}\left|  st\right|
\]
for all $\left|  s\right|  ,\left|  t\right|  <\delta$.
\end{definition}

\begin{lemma}
\label{LemmaClose}If $X$ \& $Y$ close and satisfy E1 and E2 then%
\[
d\left(  Y_{-t}X_{-t}Y_{t}X_{t}\left(  x\right)  ,x\right)  =O\left(
t^{2}\right)
\]
locally uniformly for $x\in M$.
\end{lemma}

\begin{proof}%
\begin{align*}
&  d\left(  Y_{-s}X_{-t}Y_{s}X_{t}\left(  x\right)  ,x\right) \\
&  \leq d\left(  Y_{-s}X_{-t}Y_{s}X_{t}\left(  x\right)  ,Y_{-s}X_{-t}%
X_{t}Y_{s}\left(  x\right)  \right)  +d\left(  Y_{-s}X_{-t}X_{t}Y_{s}\left(
x\right)  ,Y_{-s}Y_{s}\left(  x\right)  \right)  +d\left(  Y_{-s}Y_{s}\left(
x\right)  ,x\right) \\
&  \leq d\left(  Y_{s}X_{t}\left(  x\right)  ,X_{t}Y_{s}\left(  x\right)
\right)  \left(  1+\left|  s\right|  \Lambda_{Y}\right)  \left(  1+\left|
t\right|  \Lambda_{X}\right)  +t^{2}\Omega_{X}\left(  1+\left|  s\right|
\Lambda_{Y}\right)  +s^{2}\Omega_{Y}\\
&  \leq C_{XY}\left|  st\right|  \left(  1+\left|  s\right|  \Lambda
_{Y}\right)  \left(  1+\left|  t\right|  \Lambda_{X}\right)  +t^{2}\Omega
_{X}\left(  1+\left|  s\right|  \Lambda_{Y}\right)  +s^{2}\Omega_{Y}\leq
C\left(  \left|  st\right|  +t^{2}+s^{2}\right)
\end{align*}
where%
\[
C:=\max\left\{  C_{XY}\left(  1+\Lambda_{Y}\right)  \left(  1+\Lambda
_{X}\right)  ,\Omega_{X}\left(  1+\Lambda_{Y}\right)  ,\Omega_{Y}\right\}  .
\]
Letting $s=t$ gives the result.
\end{proof}

\begin{proposition}
If $X$ \& $Y$ close and satisfy E1 and E2 then $\left[  X,Y\right]  $ is an
arc field.
\end{proposition}

\begin{proof}
We establish the local bound on speed. The purpose of Lemma \ref{LemmaClose}
is to give $d\left(  \left[  X,Y\right]  \left(  x,t\right)  ,x\right)
=O\left(  t\right)  $ for $t\geq0$. Similarly, for $t<0$%
\begin{align*}
&  d\left(  X_{t}Y_{t}X_{-t}Y_{-t}\left(  x\right)  ,x\right) \\
&  \leq d\left(  X_{t}Y_{t}X_{-t}Y_{-t}\left(  x\right)  ,X_{t}X_{-t}\left(
x\right)  \right)  +d\left(  X_{t}X_{-t}\left(  x\right)  ,x\right) \\
&  \leq d\left(  Y_{t}X_{-t}Y_{-t}\left(  x\right)  ,X_{-t}\left(  x\right)
\right)  \left(  1+\left|  t\right|  \Lambda_{X}\right)  +t^{2}\Omega_{X}%
\end{align*}
which, using this trick again, gives%
\begin{align*}
&  \leq d\left(  X_{-t}Y_{-t}\left(  x\right)  ,Y_{-t}X_{-t}\left(  x\right)
\right)  \left(  1+\left|  t\right|  \Lambda_{X}\right)  \left(  1+\left|
t\right|  \Lambda_{Y}\right) \\
&  +t^{2}\Omega_{Y}\left(  1+\left|  t\right|  \Lambda_{Y}\right)
+t^{2}\Omega_{X}=O\left(  t^{2}\right)  \text{ since }X\text{ \& }Y\text{
close.}%
\end{align*}
Therefore
\[
d\left(  \left[  X,Y\right]  _{t}\left(  x\right)  ,x\right)  =O\left(
t\right)
\]
for both positive and negative $t$. Then since $\sqrt{\left|  t\right|  }$ is
Lipschitz except at $t=0$ we see $\left[  X,Y\right]  $ has bounded speed.
\end{proof}

\begin{example}
As in Example \ref{Banach Example} let $f,g:B\rightarrow B$ be Lipschitz
vector fields on a Banach space $B$, and let $X$ and $Y$ be their
corresponding arc fields%
\begin{align*}
X\left(  x,t\right)   &  :=x+tf\left(  x\right) \\
Y\left(  x,t\right)   &  :=x+tg\left(  x\right)
\end{align*}
It is easy to check $X$ \& $Y$ close:%
\begin{align*}
&  d\left(  Y_{s}X_{t}\left(  x\right)  ,X_{t}Y_{s}\left(  x\right)  \right)
\\
&  =\left\|  x+tf\left(  x\right)  +sg\left(  x+tf\left(  x\right)  \right)
-\left[  x+sg\left(  x\right)  +tf\left(  x+sg\left(  x\right)  \right)
\right]  \right\| \\
&  \leq\left|  t\right|  \left\|  f\left(  x\right)  -f\left(  x+sg\left(
x\right)  \right)  \right\|  +\left|  s\right|  \left\|  g\left(  x+tf\left(
x\right)  \right)  -g\left(  x\right)  \right\| \\
&  \leq\left|  t\right|  K_{f}\left\|  x-\left(  x+sg\left(  x\right)
\right)  \right\|  +\left|  s\right|  K_{g}\left\|  x+tf\left(  x\right)
-x\right\| \\
&  \leq\left|  st\right|  \left(  K_{f}\left\|  g\left(  x\right)  \right\|
+K_{g}\left\|  f\left(  x\right)  \right\|  \right)
\end{align*}
so $C_{XY}:=\left(  K_{f}\left\|  g\left(  x\right)  \right\|  +K_{g}\left\|
f\left(  x\right)  \right\|  \right)  $.

Therefore, even though the vector fields may not be smooth, so their Lie
bracket is undefined, their metric space bracket is meaningful and will give
us geometric information as we shall see in Theorem \ref{FrobeniusThm}.
\end{example}

\begin{definition}
If $X$ and $Y$ are arc fields on $M$ then define $X+Y$ to be the arc field on
$M$ given by%
\[
\left(  X+Y\right)  _{t}\left(  x\right)  :=Y_{t}X_{t}\left(  x\right)
\text{.}%
\]
For any function $a:M\rightarrow\mathbb{R}$ define the arc field $aX$ by%
\begin{equation}
aX\left(  x,t\right)  :=X\left(  x,a\left(  x\right)  t\right)  \text{.}
\label{aX}%
\end{equation}
If $a$ is Lipschitz, then $aX$ is an arc field.
\end{definition}

To be fastidiously precise we need to define $aX_{x}\left(  t\right)  $ for
all $t\in\left[  -1,1\right]  $ so technically we must specify%
\begin{equation}
aX\left(  x,t\right)  :=\left\{
\begin{array}
[c]{c}%
X\left(  x,a\left(  x\right)  t\right) \\
X\left(  x,1\right) \\
X\left(  x,-1\right) \\
x
\end{array}
\right.
\begin{array}
[c]{c}%
-\frac{1}{\left|  a\left(  x\right)  \right|  }\leq t\leq\frac{1}{\left|
a\left(  x\right)  \right|  }\\
t>1/\left|  a\left(  x\right)  \right| \\
t<-1/\left|  a\left(  x\right)  \right| \\
\text{for }-1\leq t\leq1
\end{array}%
\begin{array}
[c]{c}%
\left\}
\begin{array}
[c]{c}%
\text{ }\\
\text{when }a\left(  x\right)  \neq0\\
\text{ }%
\end{array}
\right. \\
\text{if }a\left(  x\right)  =0
\end{array}
\label{finickyArcDef}%
\end{equation}
using the trick from Example \ref{Flow<-->ArcField}. Again, we will not burden
ourselves with this detail; in all cases our concern with the properties of an
arc field $X_{x}\left(  t\right)  $ is only near $t=0$.

It is a simple definition check to prove $aX$ is an arc field when $a$ is
Lipschitz, since $aX_{x}\left(  t\right)  =X_{x}\left(  a\left(  x\right)
t\right)  $ is Lipschitz in $t$ if $X_{x}\left(  t\right)  $ is: assuming
$a\left(  x\right)  \neq0$,
\begin{align*}
\rho_{aX}\left(  x\right)   &  :=\sup_{s\neq t}\frac{d\left(  X_{x}\left(
a\left(  x\right)  s\right)  ,X_{x}\left(  a\left(  x\right)  t\right)
\right)  }{\left|  s-t\right|  }=\sup_{s\neq t}\frac{d\left(  X_{x}\left(
s\right)  ,X_{x}\left(  t\right)  \right)  }{\left|  \frac{s}{a\left(
x\right)  }-\frac{t}{a\left(  x\right)  }\right|  }\\
&  =a\left(  x\right)  \sup_{s\neq t}\frac{d\left(  X_{x}\left(  s\right)
,X_{x}\left(  t\right)  \right)  }{\left|  s-t\right|  }=a\left(  x\right)
\rho_{X}\left(  x\right)
\end{align*}
so%
\begin{align*}
\rho_{aX}\left(  x,r\right)   &  :=\sup_{y\in B\left(  x,r\right)  }\left\{
\rho_{aX}\left(  y\right)  \right\}  =\sup_{y\in B\left(  x,r\right)
}\left\{  a\left(  y\right)  \rho\left(  y\right)  \right\} \\
&  \leq\left(  a\left(  x\right)  +rK_{a}\right)  \rho_{X}\left(  x,r\right)
<\infty\text{.}%
\end{align*}

Now we have the beginnings of a linear structure associated with $M$. For
instance, expressions such as $X-Y$ make sense:%
\[
X-Y:=X+(-1)Y
\]
where $-1$ is a constant function on $M$. Further, $0$ is an arc field defined
as the constant map%
\[
0\left(  x,t\right)  :=x
\]
and satisfies $0+X=X=X+0$ for any $X$. Notice from the definition, we have
$\left[  X,Y\right]  =-\left[  Y,X\right]  $. Another trivial definition check
shows this multiplication is associative and commutative:%
\[
\left(  a\cdot b\right)  X=a\left(  bX\right)  \qquad\text{and}\qquad\left(
a\cdot b\right)  X=\left(  b\cdot a\right)  X
\]
where $\cdot$ denotes multiplication of functions.

\begin{proposition}
\label{PropX+Y E1&2}Assume $X$ \& $Y$ close and satisfy E1 and E2. Then their
sum $X+Y$ satisfies E1 and E2.
\end{proposition}

\begin{proof}
Checking Condition E1:%
\begin{align*}
&  d\left(  \left(  X+Y\right)  _{t}\left(  x\right)  ,\left(  X+Y\right)
_{t}\left(  y\right)  \right) \\
&  =d\left(  Y_{t}X_{t}\left(  x\right)  ,Y_{t}X_{t}\left(  y\right)  \right)
\leq d\left(  X_{t}\left(  x\right)  ,X_{t}\left(  y\right)  \right)  \left(
1+\left|  t\right|  \Lambda_{Y}\right) \\
&  \leq d\left(  x,y\right)  \left(  1+\left|  t\right|  \Lambda_{X}\right)
\left(  1+\left|  t\right|  \Lambda_{Y}\right)  \leq d\left(  x,y\right)
\left(  1+\left|  t\right|  \left(  \Lambda_{X}+\Lambda_{Y}\right)
+t^{2}\Lambda_{X}\Lambda_{Y}\right) \\
&  \leq d\left(  x,y\right)  \left(  1+\left|  t\right|  \Lambda_{X+Y}\right)
\end{align*}
where $\Lambda_{X+Y}:=\Lambda_{X}+\Lambda_{Y}+\Lambda_{X}\Lambda_{Y}<\infty$.

Condition E2:%
\begin{align}
&  d\left(  \left(  X+Y\right)  _{s+t}\left(  x\right)  ,\left(  X+Y\right)
_{t}\left(  X+Y\right)  _{s}\left(  x\right)  \right) \nonumber\\
&  =d\left(  Y_{s+t}X_{s+t}\left(  x\right)  ,Y_{t}X_{t}Y_{s}X_{s}\left(
x\right)  \right) \nonumber\\
&  \leq d\left(  Y_{s+t}X_{s+t}\left(  x\right)  ,Y_{t}Y_{s}X_{s+t}\left(
x\right)  \right)  +d\left(  Y_{t}Y_{s}X_{s+t}\left(  x\right)  ,Y_{t}%
X_{t}Y_{s}X_{s}\left(  x\right)  \right) \nonumber\\
&  \leq\left|  st\right|  \Omega_{Y}+d\left(  Y_{s}X_{s+t}\left(  x\right)
,X_{t}Y_{s}X_{s}\left(  x\right)  \right)  \left(  1+\left|  t\right|
\Lambda_{Y}\right) \nonumber\\
&  \leq\left|  st\right|  \Omega_{X}+\left[  d\left(  Y_{s}X_{s+t}\left(
x\right)  ,Y_{s}X_{t}X_{s}\left(  x\right)  \right)  +d\left(  Y_{s}%
X_{t}\left(  y\right)  ,X_{t}Y_{s}\left(  y\right)  \right)  \right]  \left(
1+t\Lambda_{X}\right)  \label{SumE12proof5}%
\end{align}
where $y:=X_{s}\left(  x\right)  $. Notice
\begin{align*}
d\left(  Y_{s}X_{s+t}\left(  x\right)  ,Y_{s}X_{t}X_{s}\left(  x\right)
\right)   &  \leq d\left(  X_{s+t}\left(  x\right)  ,X_{t}X_{s}\left(
x\right)  \right)  \left(  1+\left|  s\right|  \Lambda_{Y}\right) \\
&  \leq\left|  st\right|  \Omega_{X}\left(  1+\left|  s\right|  \Lambda
_{Y}\right)  =O\left(  \left|  st\right|  \right)
\end{align*}
and the last summand of $\left(  \text{\ref{SumE12proof5}}\right)  $is also
$O\left(  \left|  st\right|  \right)  $ since $X$ \& $Y$ close, so E2 is satisfied.
\end{proof}

So in this case, the flow $H$ generated by $X+Y$ is computable with Euler
curves as%
\begin{equation}
H\left(  x,t\right)  =\underset{n\rightarrow\infty}{\lim}\left(  X+Y\right)
_{t/n}^{\left(  n\right)  }\left(  x\right)  =\underset{n\rightarrow\infty
}{\lim}\left(  Y_{t/n}X_{t/n}\right)  ^{\left(  n\right)  }\left(  x\right)
\text{.} \label{X+Y Flow}%
\end{equation}
Therefore, this definition of $X+Y$ using compositions is a direct
generalization of the concept of adding vector fields on a differentiable
manifold (see \cite[Section 4.1A]{abrahamMarsden}). One of the inspirations
for this paper, \cite{Colombo} introduced the sum of semigroups on a metric
space in the same spirit as defined here, with commensurable conditions.

When $X$ \& $Y$ close and satisfy E1 and E2, we also have $\left(  X+Y\right)
\sim\left(  Y+X\right)  $ since%
\[
\left(  Y_{t/n}X_{t/n}\right)  ^{\left(  n\right)  }=Y_{t/n}\left(
X_{t/n}Y_{t/n}\right)  ^{\left(  n-1\right)  }X_{t/n}%
\]
whence both arc fields $X+Y$ and $Y+X$ are (locally uniformly) tangent to the
flow $H$ using $\left(  \ref{X+Y Flow}\right)  $.

\begin{proposition}
\label{Prop_aXE1&2}If $X$ satisfies E1 and E2 and $a:M\rightarrow\mathbb{R}$
is a Lipschitz function, then $aX$ satisfies E1 and E2.
\end{proposition}

\begin{proof}
E1:
\begin{align*}
&  d\left(  aX_{x}\left(  t\right)  ,aX_{y}\left(  t\right)  \right) \\
&  =d\left(  X_{x}\left(  a\left(  x\right)  t\right)  ,X_{y}\left(  a\left(
y\right)  t\right)  \right) \\
&  \leq d\left(  X_{x}\left(  a\left(  x\right)  t\right)  ,X_{x}\left(
a\left(  y\right)  t\right)  \right)  +d\left(  X_{x}\left(  a\left(
y\right)  t\right)  ,X_{y}\left(  a\left(  y\right)  t\right)  \right) \\
&  \leq\left|  a\left(  x\right)  t-a\left(  y\right)  t\right|  \rho\left(
x\right)  +d\left(  x,y\right)  \left(  1+a\left(  y\right)  \left|  t\right|
\Lambda_{X}\right) \\
&  \leq d\left(  x,y\right)  \left(  K_{a}\left|  t\right|  \rho\left(
x\right)  +1+a\left(  y\right)  \left|  t\right|  \Lambda_{X}\right)
=d\left(  x,y\right)  \left(  1+\left|  t\right|  \Lambda_{aX}\right)
\end{align*}
where $\Lambda_{aX}:=K_{a}\rho\left(  x\right)  +a\left(  y\right)
\Lambda_{X}<\infty.$

E2: For all $x_{0}\in M$ and $\delta>0$ we know $a$ is bounded by some $A>0$
on $B\left(  x_{0},\delta\right)  $ since $a$ is Lipschitz.%
\begin{align*}
&  d\left(  aX_{x}\left(  s+t\right)  ,aX_{aX_{x}\left(  s\right)  }\left(
t\right)  \right) \\
&  =d\left(  X_{x}\left(  a\left(  x\right)  \left(  s+t\right)  \right)
,X_{X_{x}\left(  a\left(  x\right)  s\right)  }\left(  a\left(  X_{x}\left(
a\left(  x\right)  s\right)  \right)  t\right)  \right) \\
&  \leq d\left(  X_{x}\left(  a\left(  x\right)  \left(  s+t\right)  \right)
,X_{X_{x}\left(  a\left(  x\right)  s\right)  }\left(  a\left(  x\right)
t\right)  \right) \\
&  +d\left(  X_{X_{x}\left(  a\left(  x\right)  s\right)  }\left(  a\left(
x\right)  t\right)  ,X_{X_{x}\left(  a\left(  x\right)  s\right)  }\left(
a\left(  X_{x}\left(  a\left(  x\right)  s\right)  \right)  t\right)  \right)
\\
&  \leq a\left(  x\right)  \left|  s\right|  \cdot a\left(  x\right)  \left|
t\right|  \Omega_{X}+\rho\cdot\left|  a\left(  x\right)  t-a\left(
X_{x}\left(  a\left(  x\right)  s\right)  \right)  t\right| \\
&  \leq\left|  st\right|  \left[  a\left(  x\right)  \right]  ^{2}\Omega
_{X}+\left|  t\right|  \rho K_{a}d\left(  x,X_{x}\left(  a\left(  x\right)
s\right)  \right) \\
&  \leq\left|  st\right|  \left[  a\left(  x\right)  \right]  ^{2}\Omega
_{X}+\left|  st\right|  \rho^{2}K_{a}a\left(  x\right)  \leq\left|  st\right|
\Omega_{aX}%
\end{align*}
where $\Omega_{aX}:=A^{2}\Omega_{X}+\rho^{2}K_{a}A$.
\end{proof}

Combining these results gives

\begin{theorem}
\label{ThmaX+bYflow}If $a$ and $b$ are locally Lipschitz functions and $X$ \&
$Y$ close and satisfy E1 and E2, then $aX+bY$ is an arc field which satisfies
E1 and E2 and so has a unique local flow.

If in addition $a$ and $b$ are globally Lipschitz and $X$ and $Y$ have linear
speed growth, then $aX+bY$ generates a unique flow.
\end{theorem}

\begin{proof}
We haven't proven $aX$ and $bY$ close, but this is a straightforward
definition check, as is the fact that $aX+bY$ has linear speed growth.
\end{proof}

Local flows have the following useful linearity property:

\begin{proposition}
\label{FlowLinearity}If $F$ is a local flow then interpreting $F$ as an arc
field we can perform the following operations:

1. if $a$ and $b$ are constant then\quad$aF+bF=\left(  a+b\right)  F$

2. if $a$ and $b$ are real functions then $\left(  aF+bF\right)  _{t}\left(
x\right)  =\left(  a+b\circ\left(  aF\right)  _{t}\right)  F_{t}\left(
x\right)  $.
\end{proposition}

\begin{proof}
This is another obvious definition check:%
\[%
\begin{array}
[c]{l}%
2.\quad\left(  aF+bF\right)  _{t}\left(  x\right)  =\left(  bF\right)
_{t}\left(  aF\right)  _{t}\left(  x\right)  =F_{b\left(  \left(  aF\right)
_{t}\left(  x\right)  \right)  t}F_{a\left(  x\right)  t}\left(  x\right) \\
\quad\quad\quad=F_{\left(  a\left(  x\right)  +\left(  b\circ\left(
aF\right)  _{t}\right)  \left(  x\right)  \right)  t}\left(  x\right)
=\left(  a+b\circ\left(  aF\right)  _{t}\right)  F_{t}\left(  x\right)
\end{array}
\]
and 1. follows from 2.
\end{proof}

\section{Contravariance\label{PushFwdSect}}

If $\phi:M_{1}\rightarrow M_{2}$ is a lipeomorphism (i.e., an invertible
Lipschitz map with Lipschitz inverse), then the pull-back of an arc field $X$
on $M_{2}$ is the arc field $\phi^{\ast}X$ on $M_{1}$ given by%
\[
\phi^{\ast}X\left(  x,t\right)  :=\phi^{-1}\left(  X\left(  \phi\left(
x\right)  ,t\right)  \right)
\]
or in other notation,%
\[
\left(  \phi^{\ast}X\right)  _{t}\left(  x\right)  =\phi^{-1}X_{t}\phi\left(
x\right)
\]
which is a direct analog of the pull-back of a vector field on a manifold
using curves to represent vectors. The definition for flows is identical,
replacing $X$ with $F$. The pull-back to $M_{1}$ of a solution $\sigma$ to an
arc field on $M_{2}$ is analogous:%
\[
\left(  \phi^{\ast}\sigma\right)  _{x}\left(  t\right)  :=\phi^{-1}\left(
\sigma_{\phi\left(  x\right)  }\left(  t\right)  \right)  \text{.}%
\]
The pull-back of a function $a:M_{2}\rightarrow\mathbb{R}$ is the function
$\phi^{\ast}a:M_{1}\rightarrow\mathbb{R}$ defined as $\phi^{\ast}a\left(
x\right)  :=a\left(  \phi\left(  x\right)  \right)  $.

\begin{proposition}
\label{PushFwd Sol'ns}If $\phi:M_{1}\rightarrow M_{2}$ is a lipeomorphism and
the arc field $X$ on $M_{2}$ has unique solutions then $\phi^{\ast}X$ has
unique solutions. The solutions to $\phi^{\ast}X$ are the pull-backs of
solutions to $X$.
\end{proposition}

\begin{proof}
This is obvious: if $F$ is a local flow for $X$ then%
\begin{align*}
&  d\left(  \phi^{\ast}X\left(  \phi^{\ast}F\left(  x,t\right)  ,s\right)
,\phi^{\ast}F\left(  x,t+s\right)  \right) \\
&  =d\left(  \phi^{-1}X\left[  \phi\phi^{-1}F\left(  \phi\left(  x\right)
,t\right)  ,s\right]  ,\phi^{-1}F\left(  \phi\left(  x\right)  ,t+s\right)
\right) \\
&  =d\left(  \phi^{-1}X\left[  F\left(  \phi\left(  x\right)  ,t\right)
,s\right]  ,\phi^{-1}F\left(  \phi\left(  x\right)  ,t+s\right)  \right) \\
&  \leq K_{\phi}d\left(  X\left[  F\left(  \phi\left(  x\right)  ,t\right)
,s\right]  ,F\left(  \phi\left(  x\right)  ,t+s\right)  \right)  =K_{\phi
}o\left(  s\right)  =o\left(  s\right)
\end{align*}
so $\phi^{\ast}F$ is a flow (solution) for $\phi^{\ast}X$.

Similarly if $\sigma$ is a solution to $\phi^{\ast}X$ then $\left(  \phi
^{-1}\right)  ^{\ast}\sigma$ is a solution to $X$ so by uniqueness there can
be only one such $\sigma$.
\end{proof}

The push-forward of any function, curve or flow is defined similarly, e.g.,%
\[
\phi_{\ast}F\left(  x,t\right)  :=\phi\left(  F\left(  \phi^{-1}\left(
x\right)  ,t\right)  \right)  \text{.}%
\]

It is easy to check push-forward is covariant (i.e., $\left(  \phi\circ
\psi\right)  _{\ast}=\phi_{\ast}\circ\psi_{\ast}$) and pull-back is
contravariant (i.e., $\left(  \phi\circ\psi\right)  ^{\ast}=\psi^{\ast}%
\circ\phi^{\ast}$). It is also clear that push-forward and pull-back are
inverse operations and Proposition \ref{PushFwd Sol'ns} holds \textit{mutatis
mutandis} for push-forward in place of pull-back.

\begin{proposition}
[Linearity of Pull-back]\label{PullbackLinear}If $X$ and $Y$ are arc fields on
$M$ and $\phi:M_{1}\rightarrow M_{2}$ is a lipeomorphism, then

$\left(  i\right)  $ $\phi^{\ast}\left(  X+Y\right)  =\phi^{\ast}\left(
X\right)  +\phi^{\ast}\left(  Y\right)  $

$\left(  ii\right)  $ $\phi^{\ast}\left(  aX\right)  =\left(  a\circ
\phi\right)  \phi^{\ast}\left(  X\right)  =\phi^{\ast}\left(  a\right)
\phi^{\ast}\left(  X\right)  $.
\end{proposition}

\begin{proof}
Trivial definition check.
\end{proof}

Since the pull-back and linearity are established for arc fields, we can now
explore another characterization of the bracket. In the context of $M$ being a
smooth manifold, let $F$ and $G$ be local flows generated by smooth vector
fields $f$ and $g$. There it is well known the following ``dynamic''
characterization of the Lie bracket is equivalent to the asymptotic
characterization%
\begin{equation}
\left[  f,g\right]  =\left.  \frac{d}{dt}\left(  F_{t}\right)  ^{\ast
}g\right|  _{t=0}\text{.} \label{LieDerivative}%
\end{equation}
Using%
\[
\left.  \frac{d}{dt}\left(  F_{t}\right)  ^{\ast}g\right|  _{t=0}%
=\underset{t\rightarrow0}{\lim}\frac{\left(  F_{t}\right)  ^{\ast}g-g}%
{t}=\left[  f,g\right]
\]
for inspiration, we return to the context of metric spaces where, with $F$ and
$G$ viewed as arc fields, their bracket $\left[  F,G\right]  $ is defined, and
then%
\begin{align}
F_{t}^{\ast}G_{t}\left(  x\right)   &  =\left(  t\left[  F,G\right]
+G\right)  _{t}\left(  x\right)  \text{\qquad\qquad for }t\geq0\text{ and
}\label{LieD=BracketArcFieldversion}\\
F_{s}^{\ast}G_{s}\left(  x\right)   &  =\left(  -s\left[  -F,-G\right]
-G\right)  _{-s}\left(  x\right)  \text{\qquad for }s<0
\label{LieD=BracketArcFieldversion2}%
\end{align}
which hold because%
\begin{gather*}
\left(  t\left[  F,G\right]  +G\right)  _{t}\left(  x\right)  =G_{t}\left[
F,G\right]  _{t^{2}}\left(  x\right) \\
=G_{t}G_{-t}F_{-t}G_{t}F_{t}\left(  x\right)  =F_{-t}G_{t}F_{t}\left(
x\right)  =F_{t}^{\ast}G_{t}\left(  x\right)
\end{gather*}
and%
\begin{align*}
&  \left(  -s\left[  -F,-G\right]  -G\right)  _{-s}\left(  x\right) \\
&  =G_{s}\left[  -F,-G\right]  _{s^{2}}\left(  x\right)  =G_{s}\left(
-G\right)  _{-\left|  s\right|  }\left(  -F\right)  _{-\left|  s\right|
}\left(  -G\right)  _{\left|  s\right|  }\left(  -F\right)  _{\left|
s\right|  }\left(  x\right) \\
&  =G_{s}G_{\left|  s\right|  }F_{\left|  s\right|  }G_{-\left|  s\right|
}F_{-\left|  s\right|  }\left(  x\right)  =F_{-s}G_{s}F_{s}\left(  x\right)
=F_{s}^{\ast}G_{s}\left(  x\right)  \text{.}%
\end{align*}
These facts will be used in the heart of the proof of our main result, Theorem
\ref{FrobeniusThm}, as will the following

\begin{proposition}
\label{PullbackX=X}$\left(  F_{s}\right)  ^{\ast}X\sim X$.
\end{proposition}

\begin{proof}
Using the properties of flows $F_{t}=F_{-s+t+s}=F_{-s}F_{t}F_{s}$ and
$F_{t}^{-1}=F_{t}$ we get%
\begin{align*}
&  d\left(  \left(  \left(  F_{s}\right)  ^{\ast}X\right)  _{t}\left(
x\right)  ,X_{t}\left(  x\right)  \right) \\
&  \leq d\left(  F_{-s}X_{t}F_{s}\left(  x\right)  ,F_{-s}F_{t}F_{s}\left(
x\right)  \right)  +d\left(  F_{t}\left(  x\right)  ,X_{t}\left(  x\right)
\right) \\
&  \leq e^{s\Lambda_{X}}d\left(  X_{t}\left(  y\right)  ,F_{t}\left(
y\right)  \right)  +o\left(  t\right)  =o\left(  t\right)
\end{align*}
where $y:=F_{s}\left(  x\right)  $ and the exponential comes from Theorem
\ref{ExpGrowth}.
\end{proof}

\section{Local Frobenius Theorem\label{SectionFrobThm}}

\begin{definition}
Two arc fields $X$ and $Y$ are $($locally uniformly$)$ \textbf{transverse} if
for each $x_{0}\in M$ there exists a $\delta>0$ such that%
\[
d\left(  X_{s}\left(  x\right)  ,Y_{t}\left(  x\right)  \right)  \geq
\delta\left(  \left|  s\right|  +\left|  t\right|  \right)
\]
for $\left|  t\right|  <\delta$ for all $x\in B\left(  x_{0},\delta\right)  $.
\end{definition}

\begin{example}
On the plane $\mathbb{R}^{2}$ with Euclidean norm $\left\|  \cdot\right\|  $
any two linearly independent vectors $u,v\in\mathbb{R}^{2}$ give us the
transverse arc fields%
\[
X_{t}\left(  x\right)  :=x+tu\qquad\text{and}\qquad Y_{t}\left(  x\right)
:=x+tv\text{.}%
\]
To check this, it is easiest to define a new norm on $\mathbb{R}^{2}$ by%
\[
\left\|  x\right\|  _{uv}:=\left|  x_{1}\right|  +\left|  x_{2}\right|
\]
where $x=x_{1}u+x_{2}v$ and $x_{1},x_{2}\in\mathbb{R}$. Since all norms on
$\mathbb{R}^{2}$ are metrically equivalent there must exist a constant $C>0$
such that $\left\|  x\right\|  _{uv}\leq C\left\|  x\right\|  $ for all
$x\in\mathbb{R}^{2}$. Then taking $\delta:=\frac{1}{C}$%
\[
d\left(  X_{s}\left(  x\right)  ,Y_{t}\left(  x\right)  \right)  =\left\|
su-tv\right\|  \geq\delta\left\|  su-tv\right\|  _{uv}=\delta\left(  \left|
s\right|  +\left|  t\right|  \right)  \text{.}%
\]

A localization argument shows any pair of continuous vector fields $f$ and $g$
on a differentiable manifold give transverse arc fields if $f$ and $g$ are
non-colinear at each point.
\end{example}

A (2-dimensional) \textbf{surface} is a 2-dimensional topological manifold,
i.e., locally homeomorphic to $\mathbb{R}^{2}$.

For any subset $N\subset M$ and element $x\in M$ the \textbf{distance} from
$x$ to $N$ is defined as
\[
d\left(  x,N\right)  :=\inf\left\{  d\left(  x,y\right)  :y\in N\right\}
\text{.}%
\]
This function $d$ is not a metric, obviously, but it does satisfy the triangle
inequality:%
\[
d\left(  x,N\right)  \leq d\left(  x,y\right)  +d\left(  y,N\right)
\]
for all $x,y\in M$.

\begin{definition}
A surface $S\subset M$ is an \textbf{integral surface} of two arc fields $X$
and $Y$ if for any Lipschitz functions $a,b:M\rightarrow\mathbb{R}$ then $S$
is \textbf{locally uniformly tangent} to $aX+bY$ for $x\in S$, i.e.,%
\[
d\left(  \left(  aX+bY\right)  _{t}\left(  x\right)  ,S\right)  =o\left(
t\right)
\]
locally uniformly in $x$. Locally uniform tangency is denoted $S\sim aX+bY$.
\end{definition}

\begin{theorem}
\label{FrobeniusThm}Assume $X$ \& $Y$ close, are transverse, and satisfy E1
and E2 on a locally complete metric space $M$. Let $F$ and $G$ be the local
flows of $X$ and $Y$. If $\left[  F,G\right]  \sim aX+bY$ $($locally uniform
tangency$)$ for some Lipschitz functions $a,b:M\rightarrow\mathbb{R}$, then
for each $x_{0}\in M$ there exists an integral surface $S$ through $x_{0}$.
\end{theorem}

\begin{proof}
It may be beneficial to review the outline of this proof from the third
paragraph of the introduction. The metric space constructs of the previous
sections will now be inserted into the manifold outline. A rigorous
verification of the analytic estimates requires some tedious, but
straightforward, calculations detailed here.

Define
\[
S:=\left\{  F_{t}G_{s}\left(  x_{0}\right)  :\left|  s\right|  ,\left|
t\right|  <\delta\right\}
\]
where $\delta>0$ is chosen small enough for $S$ to be a well-defined surface
(Figure \ref{FigFrobProof1}).
%TCIMACRO{\FRAME{ftbphFU}{3.6019in}{2.3722in}{0pt}{\Qcb{integral surface $S$}%
%}{\Qlb{FigFrobProof1}}{frobproof1.eps}{\special{ language "Scientific Word";
%type "GRAPHIC";  maintain-aspect-ratio TRUE;  display "USEDEF";
%valid_file "F";  width 3.6019in;  height 2.3722in;  depth 0pt;
%original-width 5.1923in;  original-height 3.4117in;  cropleft "0";
%filename 'FrobProof1.eps';file-properties "XNPEU";}} }%
%BeginExpansion
\begin{figure}
[ptbh]
\begin{center}
\includegraphics[
height=2.3722in,
width=3.6019in
]%
{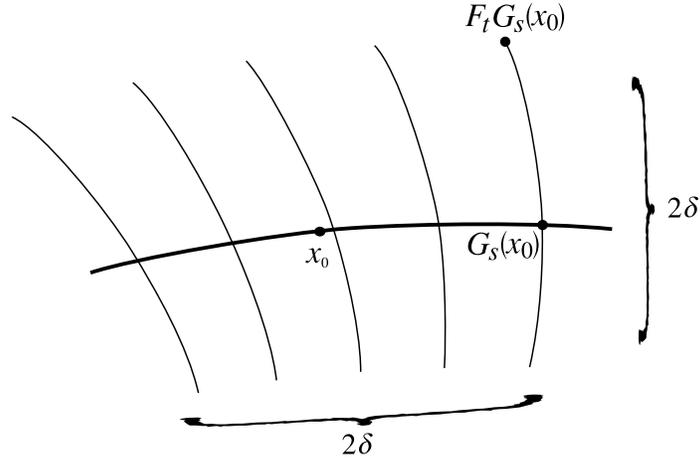}%
\caption{integral surface $S$}%
\label{FigFrobProof1}%
\end{center}
\end{figure}
%EndExpansion
I.e., $F_{t_{1}}G_{s_{1}}\left(  x_{0}\right)  =F_{t_{2}}G_{s_{2}}\left(
x_{0}\right)  $ implies $t_{1}=t_{2}$ and $s_{1}=s_{2}$ so%
\[
\phi:\left(  -\delta,\delta\right)  \times\left(  -\delta,\delta\right)
\subset\mathbb{R}^{2}\rightarrow S\subset M
\]
defined by $\phi\left(  s,t\right)  :=F_{t}G_{s}\left(  x_{0}\right)  $ is a
homeomorphism. Finding such a $\delta$ is possible since $X$ and $Y$ are
transverse. To see this, assume the contrary. Then there are different choices
of $s_{i}$ and $t_{i}$ which give $F_{t_{1}}G_{s_{1}}\left(  x_{0}\right)
=F_{t_{2}}G_{s_{2}}\left(  x_{0}\right)  $ which implies $G_{s_{1}}\left(
x_{0}\right)  =F_{t_{3}}G_{s_{2}}\left(  x_{0}\right)  $ and letting
$y:=G_{s_{2}}\left(  x_{0}\right)  $ we must also then have
\begin{equation}
F_{t}\left(  y\right)  =G_{s}\left(  y\right)  \text{.}\label{FrobProof3}%
\end{equation}
If this contrary assumption were true, then for all $\varepsilon>0$ there
would exist $s$ and $t$ with $\left|  s\right|  ,\left|  t\right|
<\varepsilon$ such that $\left(  \ref{FrobProof3}\right)  $ holds. Since $X$
and $Y$ are transverse, this cannot be so.

We will show $S$ is the desired integral surface through $x_{0}$. Assume
$\delta$ is also chosen small enough so throughout $S$ the functions $\left|
a\right|  $ and $\left|  b\right|  $ are bounded, while the constants
$\Lambda$, $\Omega$, and $\rho$ for $X$ and $Y$ hold uniformly, and that the
closure of $B\left(  x,2\delta\left(  \rho+1\right)  \right)  $ is complete.
This is possible because $F$ and $G$ have locally bounded speeds, since $X$
and $Y$ do.

Notice $S\sim X$ by construction, but it is not immediately clear $S\sim
a^{\prime}X+b^{\prime}Y$ for arbitrarily chosen $a^{\prime},b^{\prime}%
\in\mathbb{R}$. Notice we can use
\[
a^{\prime}X+b^{\prime}Y\sim a^{\prime}F+b^{\prime}G\sim b^{\prime}G+a^{\prime
}F\sim b^{\prime}Y+a^{\prime}X
\]
and so we will show $S\sim a^{\prime}F+b^{\prime}G$. We need to show this is
true for an arbitrary point $z\in S,$ so assume $z:=F_{t}G_{s}\left(
x_{0}\right)  $ for some $s$ and $t\in\mathbb{R}$. Notice by the construction
of $S$ we have $S\sim a^{\prime\prime}F+b^{\prime\prime}G$ at $x:=G_{s}\left(
x_{0}\right)  $ for an arbitrary choice of Lipschitz functions $a^{\prime
\prime}$ and $b^{\prime\prime}$ since $a^{\prime\prime}F+b^{\prime\prime}G\sim
b^{\prime\prime}G+a^{\prime\prime}F$ and%
\begin{align*}
&  \left(  b^{\prime\prime}G+a^{\prime\prime}F\right)  _{h}\left(  x\right) \\
&  =F_{a^{\prime\prime}\left(  G_{b^{\prime\prime}\left(  x\right)  h}\left(
x\right)  \right)  h}G_{b^{\prime\prime}\left(  x\right)  h}\left(  x\right)
=F_{a^{\prime\prime}\left(  G_{b^{\prime\prime}\left(  x\right)  h}\left(
x\right)  \right)  h}G_{b^{\prime\prime}\left(  x\right)  h}\left(  x\right)
\\
&  =F_{a^{\prime\prime}\left(  G_{b^{\prime\prime}\left(  x\right)  h}\left(
x\right)  \right)  h}G_{b^{\prime\prime}\left(  x\right)  h}G_{s}\left(
x_{0}\right)  \in S
\end{align*}
when $h$ is small.

($x_{0},x,z,s$ and $t$ are now fixed for the remainder of the proof; however,
we only explicitly check the case $t>0$, indicating the changes where needed
to check the $t<0$ case.)

If we prove%
\begin{equation}
\left(  F_{t}\right)  ^{\ast}\left(  a^{\prime}F+b^{\prime}G\right)  \sim
S\text{\qquad at\qquad}x=G_{s}\left(  x_{0}\right)  \label{FrobProof4}%
\end{equation}
this will prove $S\sim a^{\prime}F+b^{\prime}G$ at $z$, since the push-forward
$\left(  F_{t}\right)  _{\ast}$ and the pull-back $\left(  F_{t}\right)
^{\ast}$ are inverse and local lipeomorphisms and so preserve tangency. See
Figure \ref{FigFrobProof2}.%
%TCIMACRO{\FRAME{ftbphFU}{3.9392in}{2.1689in}{0pt}{\Qcb{pull-back to
%$G_{s}\left(  x_{0}\right)  $}}{\Qlb{FigFrobProof2}}{frobprooffig2.eps}%
%width 3.9392in;  height 2.1689in;  depth 0pt;  original-width 5.4059in;
%original-height 2.9654in;  cropleft "0";  croptop "1";  cropright "1";
%cropbottom "0";  filename 'FrobProofFig2.eps';file-properties "XNPEU";}} }%
%BeginExpansion
\begin{figure}
[ptbh]
\begin{center}
\includegraphics[
height=2.1689in,
width=3.9392in
]%
{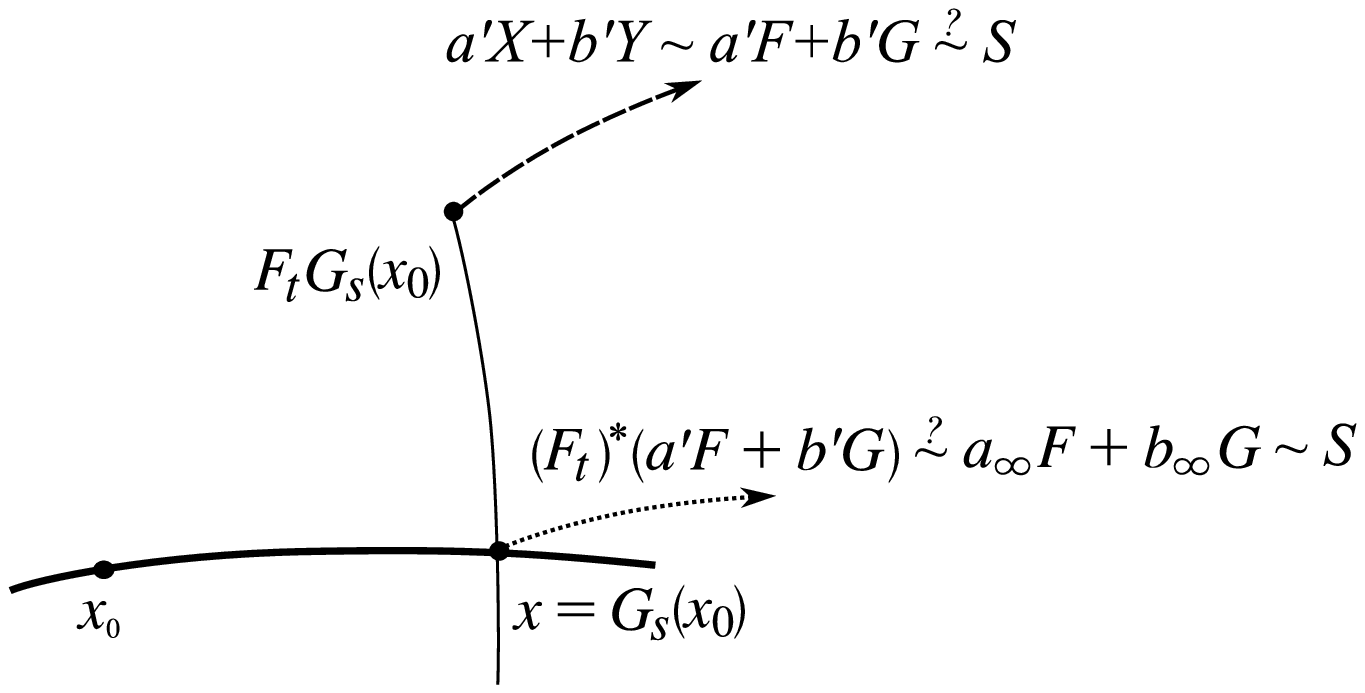}%
\caption{pull-back to $G_{s}\left(  x_{0}\right)  $}%
\label{FigFrobProof2}%
\end{center}
\end{figure}
%EndExpansion

Restating $\left(  \text{\ref{LieD=BracketArcFieldversion}}\right)  $:%
\[
F_{t}^{\ast}G_{t}\left(  x\right)  =\left(  t\left[  F,G\right]  +G\right)
_{t}\left(  x\right)
\]
so%
\begin{equation}
F_{t/n}^{\ast}G_{t/n}\left(  x\right)  =\left(  \tfrac{t}{n}\left[
F,G\right]  +G\right)  _{t/n}\left(  x\right)  \label{FrobProof20}%
\end{equation}
for our previously fixed small $t\geq0$ and arbitrary positive integer
$n\in\mathbb{N}$. (For $t<0$ use $\left(
\text{\ref{LieD=BracketArcFieldversion2}}\right)  $ instead.) For any arc
fields $Z$ and $\overline{Z}$ clearly%
\begin{gather}
d\left(  Z_{t}\left(  x\right)  ,\overline{Z}_{t}\left(  x\right)  \right)
=o\left(  t\right)  \qquad\text{implies}\nonumber\\
d\left(  \left(  tZ\right)  _{t}\left(  x\right)  ,\left(  t\overline
{Z}\right)  _{t}\left(  x\right)  \right)  =d\left(  \left(  Z\right)
_{t^{2}}\left(  x\right)  ,\left(  \overline{Z}\right)  _{t^{2}}\left(
x\right)  \right)  =o\left(  t^{2}\right)  \label{FrobProof30}%
\end{gather}
and so%
\begin{gather}
\left[  F,G\right]  \sim aF+bG\qquad\text{implies}\nonumber\\
d\left(  \left(  \tfrac{t}{n}\left[  F,G\right]  \right)  _{t/n}\left(
x\right)  ,\left(  \left(  \tfrac{t}{n}\left(  aF+bG\right)  \right)  \right)
_{t/n}\left(  x\right)  \right)  =o\left(  \tfrac{1}{n^{2}}\right)
\label{FrobProof40}%
\end{gather}
since $t$ is fixed.

We use these facts to establish $\left(  \text{\ref{FrobProof4}}\right)  $,
first checking%
\[
d\left(  \left(  F_{t}^{\ast}\left(  a^{\prime}F+b^{\prime}G\right)  \right)
_{t/n}\left(  x\right)  ,S\right)  =o\left(  \tfrac{1}{n}\right)
\]
as $n\rightarrow\infty$. At the end of the proof we will replace $t/n$ by
arbitrary $r\rightarrow0$. Using the linearity of pull-back (Proposition
\ref{PullbackLinear}) we get
\begin{align}
&  d\left(  \left(  F_{t}^{\ast}\left(  a^{\prime}F+b^{\prime}G\right)
\right)  _{t/n}\left(  x\right)  ,S\right) \nonumber\\
&  =d\left(  \left(  \left(  a^{\prime}\circ F_{t}\right)  F_{t}^{\ast}\left(
F\right)  +\left(  b^{\prime}\circ F_{t}\right)  F_{t/n}^{\ast\left(
n\right)  }\left(  G\right)  \right)  _{t/n}\left(  x\right)  ,S\right)
\nonumber\\
&  =d\left(  \left(  a_{0}F+b_{0}F_{t/n}^{\ast\left(  n\right)  }\left(
G\right)  \right)  _{t/n}\left(  x\right)  ,S\right) \nonumber
\end{align}
where $a_{0}:=a^{\prime}\circ F_{t}$ and $b_{0}:=b^{\prime}\circ F_{t}$. Using
$\left(  \text{\ref{FrobProof20}}\right)  $ means this last estimate is%
\begin{align}
&  =d\left(  \left(  a_{0}F+b_{0}F_{t/n}^{\ast\left(  n-1\right)  }\left(
\tfrac{t}{n}\left[  F,G\right]  +G\right)  \right)  _{t/n}\left(  x\right)
,S\right) \nonumber\\
&  \leq d\left(  \left(  a_{0}F+b_{0}F_{t/n}^{\ast\left(  n-1\right)  }\left(
\tfrac{t}{n}\left[  F,G\right]  +G\right)  \right)  _{t/n}\left(  x\right)
,\left(  a_{0}F+b_{0}F_{t/n}^{\ast\left(  n-1\right)  }\left(  \tfrac{t}%
{n}\left(  aF+bG\right)  +G\right)  \right)  _{t/n}\left(  x\right)  \right)
\nonumber\\
&  +d\left(  \left(  a_{0}F+b_{0}F_{t/n}^{\ast\left(  n-1\right)  }\left(
\tfrac{t}{n}\left(  aF+bG\right)  +G\right)  \right)  _{t/n}\left(  x\right)
,S\right)  \text{.} \label{FrobProof45}%
\end{align}
We estimate the first term as%
\begin{align*}
&  d\left(  \left(  a_{0}F+b_{0}F_{t/n}^{\ast\left(  n-1\right)  }\left(
\tfrac{t}{n}\left[  F,G\right]  +G\right)  \right)  _{t/n}\left(  x\right)
,\left(  a_{0}F+b_{0}F_{t/n}^{\ast\left(  n-1\right)  }\left(  \tfrac{t}%
{n}\left(  aF+bG\right)  +G\right)  \right)  _{t/n}\left(  x\right)  \right)
\\
&  =d\left(  \left(  b_{0}F_{\left(  n-1\right)  t/n}^{\ast}\left(  \tfrac
{t}{n}\left[  F,G\right]  +G\right)  \right)  _{t/n}\left(  y\right)  ,\left(
b_{0}F_{\left(  n-1\right)  t/n}^{\ast}\left(  \tfrac{t}{n}\left(
aF+bG\right)  +G\right)  \right)  _{t/n}\left(  y\right)  \right)
\end{align*}
where $y:=a_{0}F_{t/n}\left(  x\right)  $%
\begin{align}
&  =d\left(  \left(  F_{\left(  n-1\right)  t/n}^{\ast}\left(  \tfrac{t}%
{n}\left[  F,G\right]  +G\right)  \right)  _{b_{0}\left(  y\right)
t/n}\left(  y\right)  ,\left(  F_{\left(  n-1\right)  t/n}^{\ast}\left(
\tfrac{t}{n}\left(  aF+bG\right)  +G\right)  \right)  _{b_{0}\left(  y\right)
t/n}\left(  y\right)  \right) \nonumber\\
&  =d\left(
\begin{array}
[c]{c}%
\left(  F_{-\left(  n-1\right)  t/n}\left(  \tfrac{t}{n}\left[  F,G\right]
+G\right)  \right)  _{b_{0}\left(  y\right)  t/n}\left(  F_{\left(
n-1\right)  t/n}\left(  y\right)  \right) \\
,\left(  F_{-\left(  n-1\right)  t/n}\left(  \tfrac{t}{n}\left(  aF+bG\right)
+G\right)  \right)  _{b_{0}\left(  y\right)  t/n}\left(  F_{\left(
n-1\right)  t/n}\left(  y\right)  \right)
\end{array}
\right) \nonumber\\
&  =d\left(
\begin{array}
[c]{c}%
\left(  F_{-\left(  n-1\right)  t/n}\left(  \tfrac{t}{n}\left[  F,G\right]
+G\right)  \right)  _{b_{0}\left(  y\right)  t/n}\left(  z\right) \\
,\left(  F_{-\left(  n-1\right)  t/n}\left(  \tfrac{t}{n}\left(  aF+bG\right)
+G\right)  \right)  _{b_{0}\left(  y\right)  t/n}\left(  z\right)
\end{array}
\right)  \label{FrobProof60}%
\end{align}
where $z:=F_{\left(  n-1\right)  t/n}\left(  y\right)  $. Then by Theorem
\ref{ExpGrowth}, $\left(  \ref{FrobProof60}\right)  $ is%
\begin{align}
&  \leq d\left(  \left(  \tfrac{t}{n}\left[  F,G\right]  +G\right)
_{b_{0}\left(  y\right)  t/n}\left(  z\right)  ,\left(  \tfrac{t}{n}\left(
aF+bG\right)  +G\right)  _{b_{0}\left(  y\right)  t/n}\left(  z\right)
\right)  e^{\Lambda_{X}\left(  n-1\right)  t/n}\nonumber\\
&  =d\left(  G_{b_{0}\left(  y\right)  t/n}\left(  \tfrac{t}{n}\left[
F,G\right]  \right)  _{b_{0}\left(  y\right)  t/n}\left(  z\right)
,G_{b_{0}\left(  y\right)  t/n}\left(  \tfrac{t}{n}\left(  aF+bG\right)
\right)  _{b_{0}\left(  y\right)  t/n}\left(  z\right)  \right)
e^{\Lambda_{X}\left(  n-1\right)  t/n}\nonumber\\
&  \leq d\left(  \left(  \tfrac{t}{n}\left[  F,G\right]  \right)
_{b_{0}\left(  y\right)  t/n}\left(  z\right)  ,\left(  \tfrac{t}{n}\left(
aF+bG\right)  \right)  _{b_{0}\left(  y\right)  t/n}\left(  z\right)  \right)
e^{\Lambda_{X}\left(  n-1\right)  t/n}e^{\Lambda_{Y}b_{0}\left(  y\right)
t/n}\nonumber\\
&  \leq r\left(  b_{0}\left(  y\right)  \left(  \tfrac{t}{n}\right)
^{2}\right)  e^{\Lambda_{X}\left(  n-1\right)  t/n+\Lambda_{Y}b_{0}\left(
y\right)  t/n}=:o_{1}\left(  \tfrac{1}{n^{2}}\right)  \label{FrobProof65}%
\end{align}
where we define%
\[
r\left(  s\right)  :=d\left(  \left[  F,G\right]  _{s}\left(  z\right)
,\left(  aF+bG\right)  _{s}\left(  z\right)  \right)  \text{.}%
\]
By the main assumption of the theorem, $r\left(  s\right)  =o\left(  s\right)
$ so notice we have $o_{1}\left(  \tfrac{1}{n^{2}}\right)  =o\left(  \tfrac
{1}{n^{2}}\right)  $ but we need to keep a careful record of this estimate as
we will be summing $n$ terms like it; the subscript distinguishes $o_{1}$ as a
specific function.

Substituting $\left(  \ref{FrobProof65}\right)  $ into $\left(
\ref{FrobProof45}\right)  $ gives%
\begin{align}
&  d\left(  \left(  F_{t}^{\ast}\left(  a^{\prime}F+b^{\prime}G\right)
\right)  _{t/n}\left(  x\right)  ,S\right) \nonumber\\
&  =d\left(  \left(  a_{0}F+b_{0}F_{t/n}^{\ast\left(  n\right)  }G\right)
_{t/n}\left(  x\right)  ,S\right) \label{FrobProof70}\\
&  \leq d\left(  \left(  a_{0}F+b_{0}F_{t/n}^{\ast\left(  n-1\right)  }\left(
\tfrac{t}{n}\left(  aF+bG\right)  +G\right)  \right)  _{t/n}\left(  x\right)
,S\right)  +o_{1}\left(  \tfrac{1}{n^{2}}\right) \nonumber\\
&  =d\left(  \left(
\begin{array}
[c]{c}%
a_{0}F+b_{0}\tfrac{t}{n}\left(  a\circ F_{\left(  n-1\right)  t/n}\right)  F\\
+b_{0}\cdot\left(  \tfrac{t}{n}\left(  b\circ F_{\left(  n-1\right)
t/n}\right)  +1\right)  F_{t/n}^{\ast\left(  n-1\right)  }G
\end{array}
\right)  _{t/n}\left(  x\right)  ,S\right)  +o_{1}\left(  \tfrac{1}{n^{2}%
}\right) \nonumber\\
&  =d\left(  \left(
\begin{array}
[c]{c}%
\left[  a_{0}+\left(  b_{0}\tfrac{t}{n}\left(  a\circ F_{\left(  n-1\right)
t/n}\right)  \right)  \circ\left(  a_{0}F_{t/n}\right)  \right]  F\\
+b_{0}\cdot\left(  \tfrac{t}{n}\left(  b\circ F_{\left(  n-1\right)
t/n}\right)  +1\right)  F_{t/n}^{\ast\left(  n-1\right)  }G
\end{array}
\right)  _{t/n}\left(  x\right)  ,S\right)  +o_{1}\left(  \tfrac{1}{n^{2}%
}\right) \nonumber\\
&  =d\left(  \left(  a_{1}F+b_{1}F_{t/n}^{\ast\left(  n-1\right)  }G\right)
_{t/n}\left(  x\right)  ,S\right)  +o_{1}\left(  \tfrac{1}{n^{2}}\right)
\label{FrobProof80}%
\end{align}
where%
\begin{align*}
a_{1}  &  :=a_{0}+\left(  b_{0}\tfrac{t}{n}\left(  a\circ F_{\left(
n-1\right)  t/n}\right)  \right)  \circ\left(  a_{0}F_{t/n}\right)
\qquad\text{and}\\
b_{1}  &  :=b_{0}\cdot\left(  \tfrac{t}{n}\left(  b\circ F_{\left(
n-1\right)  t/n}\right)  +1\right)  \text{.}%
\end{align*}
This painful calculation from the third line to the fourth line employs the
linearity of pull-back (Proposition \ref{PullbackLinear}); while the fifth
line is due to the linearity of $F$ (Proposition \ref{FlowLinearity}).

After toiling through these many complicated estimates we can relax a bit,
since the rest of the proof follows more mechanically by iterating the result
of lines $\left(  \ref{FrobProof70}\right)  $ and $\left(  \ref{FrobProof80}%
\right)  $:%
\begin{align}
&  d\left(  \left(  a_{0}F+b_{0}F_{t/n}^{\ast\left(  n\right)  }G\right)
_{t/n}\left(  x\right)  ,S\right) \nonumber\\
&  \leq d\left(  \left(  a_{1}F+b_{1}F_{t/n}^{\ast\left(  n-1\right)
}G\right)  _{t/n}\left(  x\right)  ,S\right)  +o_{1}\left(  \tfrac{1}{n^{2}%
}\right) \nonumber\\
&  \leq d\left(  \left(  a_{2}F+b_{2}F_{t/n}^{\ast\left(  n-2\right)
}G\right)  _{t/n}\left(  x\right)  ,S\right)  +o_{1}\left(  \tfrac{1}{n^{2}%
}\right)  +o_{2}\left(  \tfrac{1}{n^{2}}\right) \nonumber\\
&  \leq...\leq d\left(  \left(  a_{n}F+b_{n}G\right)  _{t/n}\left(  x\right)
,S\right)  +\overset{n}{\underset{i=1}{%
%TCIMACRO{\tsum }%
%BeginExpansion
{\textstyle\sum}
%EndExpansion
}}o_{i}\left(  \tfrac{1}{n^{2}}\right)  \label{FrobProof90}%
\end{align}
where
\begin{align*}
a_{2}  &  :=a_{1}+\left(  b_{1}\tfrac{t}{n}\left(  a\circ F_{\left(
n-2\right)  t/n}\right)  \right)  \circ\left(  a_{1}F_{t/n}\right) \\
b_{2}  &  :=b_{1}\cdot\left(  \tfrac{t}{n}\left(  b\circ F_{\left(
n-2\right)  t/n}\right)  +1\right)  \text{\qquad\qquad and in general}\\
a_{i}  &  :=a_{i-1}+\left(  b_{i-1}\tfrac{t}{n}\left(  a\circ F_{\left(
n-i\right)  t/n}\right)  \right)  \circ\left(  a_{i-1}F_{t/n}\right) \\
b_{i}  &  :=b_{i-1}\cdot\left(  \tfrac{t}{n}\left(  b\circ F_{\left(
n-i\right)  t/n}\right)  +1\right)
\end{align*}
In the region of interest the $\left|  a\right|  $ and $\left|  a_{0}\right|
$ are bounded by some $A\in\mathbb{R}$ and $\left|  b\right|  $ and $\left|
b_{0}\right|  $ are bounded by some $B\in\mathbb{R}$ so%
\begin{align*}
\left|  b_{1}\right|   &  =\left|  b_{0}\cdot\left(  \tfrac{t}{n}\left(
b\circ F_{\left(  n-1\right)  t/n}\right)  +1\right)  \right|  \leq B\left(
\tfrac{t}{n}B+1\right) \\
\left|  b_{2}\right|   &  =\left|  b_{1}\cdot\left(  \tfrac{t}{n}\left(
b\circ F_{\left(  n-1\right)  t/n}\right)  +1\right)  \right|  \leq B\left(
\tfrac{t}{n}B+1\right)  ^{2}\\
\left|  b_{i}\right|   &  \leq B\left(  \tfrac{t}{n}B+1\right)  ^{i}%
\text{\qquad and}%
\end{align*}%
\begin{align*}
\left|  a_{1}\right|   &  =\left|  a_{0}+b_{0}\tfrac{t}{n}\left(  a\circ
F_{\left(  n-1\right)  t/n}\right)  \right|  \leq A+B\tfrac{t}{n}A\\
\left|  a_{2}\right|   &  =\left|  a_{1}+b_{1}\tfrac{t}{n}\left(  a\circ
F_{\left(  n-2\right)  t/n}\right)  \right|  \leq\left(  A+B\tfrac{t}%
{n}A\right)  +B\left(  \tfrac{t}{n}B+1\right)  \tfrac{t}{n}A\\
\left|  a_{3}\right|   &  =\left|  a_{2}+b_{2}\tfrac{t}{n}\left(  a\circ
F_{\left(  n-3\right)  t/n}\right)  \right| \\
&  \leq A+B\tfrac{t}{n}A+B\left(  \tfrac{t}{n}B+1\right)  \tfrac{t}%
{n}A+B\left(  \tfrac{t}{n}B+1\right)  ^{2}\tfrac{t}{n}A\\
\left|  a_{i}\right|   &  \leq A+\tfrac{t}{n}AB\overset{i-1}{\underset{k=0}{%
%TCIMACRO{\tsum }%
%BeginExpansion
{\textstyle\sum}
%EndExpansion
}}\left(  \tfrac{t}{n}B+1\right)  ^{k}=A+\tfrac{t}{n}AB\frac{\left(  \tfrac
{t}{n}B+1\right)  ^{i}-1}{\tfrac{t}{n}B}\\
&  =A\left(  \tfrac{t}{n}B+1\right)  ^{i}\text{.}%
\end{align*}
Therefore%
\begin{align*}
\left|  b_{n}\right|   &  \leq B\left(  \tfrac{t}{n}B+1\right)  ^{n}\leq
Be^{tB}\text{\qquad and}\\
\left|  a_{n}\right|   &  \leq A\left(  \tfrac{t}{n}B+1\right)  ^{n}\leq
Ae^{tB}\text{.}%
\end{align*}

Penultimately, we need to estimate the $o_{i}\left(  \tfrac{1}{n^{2}}\right)
$. Remember from line $\left(  \ref{FrobProof65}\right)  $%
\[
o_{1}\left(  \tfrac{1}{n^{2}}\right)  :=r\left(  b_{0}\left(  y\right)
\left(  \tfrac{t}{n}\right)  ^{2}\right)  e^{\Lambda_{X}\left(  n-1\right)
t/n+\Lambda_{Y}b_{0}\left(  y\right)  t/n}%
\]
where $r\left(  s\right)  =o\left(  s\right)  $, so%
\begin{align*}
o_{2}\left(  \tfrac{1}{n^{2}}\right)   &  =r\left(  b_{1}\left(  y\right)
\left(  \tfrac{t}{n}\right)  ^{2}\right)  e^{\Lambda_{X}\left(  n-2\right)
t/n+\Lambda_{Y}b_{1}\left(  y\right)  t/n}\\
&  \leq B\left(  \tfrac{t}{n}B+1\right)  o\left(  \left(  \tfrac{t}{n}\right)
^{2}\right)  e^{\Lambda_{X}\left(  n-2\right)  t/n+\Lambda_{Y}B\left(
\tfrac{t}{n}B+1\right)  t/n}\\
o_{i}\left(  \tfrac{1}{n^{2}}\right)   &  =r\left(  b_{i-1}\left(  y\right)
\left(  \tfrac{t}{n}\right)  ^{2}\right)  e^{\Lambda_{X}\left(  n-i\right)
t/n+\Lambda_{Y}b_{i-1}\left(  y\right)  t/n}\text{.}%
\end{align*}
Consequently%
\begin{align*}
\overset{n}{\underset{i=1}{%
%TCIMACRO{\tsum }%
%BeginExpansion
{\textstyle\sum}
%EndExpansion
}}o_{i}\left(  \tfrac{1}{n^{2}}\right)   &  \leq\overset{n}{\underset{i=1}{%
%TCIMACRO{\tsum }%
%BeginExpansion
{\textstyle\sum}
%EndExpansion
}}r\left(  b_{i-1}\left(  y\right)  \left(  \tfrac{t}{n}\right)  ^{2}\right)
e^{\Lambda_{X}\left(  n-i\right)  t/n+\Lambda_{Y}B\left(  \tfrac{t}%
{n}B+1\right)  ^{i-1}t/n}\\
&  \leq o\left(  \left(  \tfrac{t}{n}\right)  ^{2}\right)  Be^{tB}\overset
{n}{\underset{i=1}{%
%TCIMACRO{\tsum }%
%BeginExpansion
{\textstyle\sum}
%EndExpansion
}}e^{\Lambda_{X}\left(  n-i\right)  t/n+\Lambda_{Y}B\left(  \tfrac{t}%
{n}B+1\right)  ^{i-1}t/n}%
\end{align*}
since $r\left(  b_{i-1}\left(  y\right)  \left(  \tfrac{t}{n}\right)
^{2}\right)  =o\left(  \left(  \tfrac{t}{n}\right)  ^{2}\right)  Be^{tB}$ for
all $i$. Therefore%
\[
\overset{n}{\underset{i=1}{%
%TCIMACRO{\tsum }%
%BeginExpansion
{\textstyle\sum}
%EndExpansion
}}o_{i}\left(  \tfrac{1}{n^{2}}\right)  \leq o\left(  \left(  \tfrac{t}%
{n}\right)  ^{2}\right)  Be^{tB}ne^{\Lambda_{X}t+\Lambda_{Y}Be^{tB}%
t/n}=o\left(  \tfrac{1}{n}\right)
\]
as $n\rightarrow\infty$. Putting this into $\left(  \ref{FrobProof90}\right)
$ gives%
\[
d\left(  \left(  F_{t}^{\ast}\left(  a^{\prime}F+b^{\prime}G\right)  \right)
_{t/n}\left(  x\right)  ,S\right)  \leq d\left(  \left(  a_{n}F+b_{n}G\right)
_{t/n}\left(  x\right)  ,S\right)  +o\left(  \tfrac{1}{n}\right)  =o\left(
\tfrac{1}{n}\right)
\]
because of the uniform bound on $\left|  a_{n}\right|  $ and $\left|
b_{n}\right|  $. To see this notice%
\[
d\left(  \left(  a_{\ast}F+b_{\ast}G\right)  _{t/n}\left(  x\right)
,S\right)  =o\left(  \tfrac{1}{n}\right)
\]
uniformly for bounded $a_{\ast}$ and $b_{\ast}$ since $a_{\ast}F+b_{\ast}G\sim
b_{\ast}G+a_{\ast}F$ and as before $\left(  b_{\ast}G+a_{\ast}F\right)
_{t}\left(  x\right)  \in S$ using the uniform $\Lambda$ and $\Omega$ derived
in the proofs of Propositions \ref{PropX+Y E1&2} and \ref{Prop_aXE1&2} (cf.
Remark \ref{RemUniformSolutions}).

Finally we need to check%
\[
d\left(  \left(  F_{t}^{\ast}\left(  a^{\prime}F+b^{\prime}G\right)  \right)
_{r}\left(  x\right)  ,S\right)  =o\left(  r\right)
\]
when $r$ is not necessarily $t/n$. We may assume $0<t<1$ and $0<r<t$ so that
$t=nr+\varepsilon$ for some $0\leq\varepsilon<r$ and integer $n$ with
$\frac{t}{r}-1<n\leq\frac{t}{r}$. Therefore the above calculations give%
\begin{align*}
&  d\left(  \left(  F_{t}^{\ast}\left(  a^{\prime}F+b^{\prime}G\right)
\right)  _{r}\left(  x\right)  ,S\right)  =d\left(  \left(  F_{\varepsilon
}^{\ast}F_{r}^{\ast\left(  n\right)  }\left(  cF+dG\right)  \right)
_{r}\left(  x\right)  ,S\right) \\
&  \leq d\left(  F_{\varepsilon}^{\ast}\left(  a_{n}F+b_{n}G\right)
_{r}\left(  x\right)  ,S\right)  +o\left(  r\right)  =o\left(  r\right)
\text{.}%
\end{align*}
\end{proof}

The $n$-dimensional corollary of this $2$-dimensional version is given in the
next section.

\begin{remark}
\label{Rem2ndObracket}In the assumptions of Theorem \ref{FrobeniusThm}
$\left[  F,G\right]  $ can be replaced with $\left[  X,Y\right]  $ when they
are tangent. Since the brackets use $\sqrt{t}$ we have $\left[  F,G\right]
\sim\left[  X,Y\right]  $ when $X$ and $Y$ are $2$nd-order tangent to their
flows, i.e.,
\begin{align*}
d\left(  X_{t}\left(  x\right)  ,F_{t}\left(  x\right)  \right)   &  =O\left(
t^{2}\right)  \qquad\qquad\text{and}\\
d\left(  Y_{t}\left(  x\right)  ,G_{t}\left(  x\right)  \right)   &  =O\left(
t^{2}\right)
\end{align*}
locally uniformly. We denote $2$nd-order local uniform tangency by $X\approx
F$. This holds, for example, when $X$ comes from a twice continuously
differentiable vector field by Taylor's theorem. But in formulating our
theorem for the nonsmooth case, the two brackets are not interchangeable.
Beware: 2nd-order tangency is ``big oh'' of $t^{2}$, not ``little oh''.

We might have chosen to define the bracket $\left[  X,Y\right]  $ using the
flows instead of the arc fields to simplify the statements of Theorem
\ref{FrobeniusThm} and those below. However it is often easier to calculate
the bracket and to check closure using arc fields instead of the flows.

In light of this remark, Theorem \ref{FrobeniusThm} gives
\end{remark}

\begin{corollary}
Assume $X$ \& $Y$ close, are transverse, and satisfy E1 and E2 on a locally
complete metric space $M$. Further assume $X$ and $Y$ are 2nd-order tangent to
their local flows $F$ and $G$. If $\left[  X,Y\right]  \sim aX+bY$ for some
Lipschitz functions $a,b:M\rightarrow\mathbb{R}$, then for each $x_{0}\in M$
there exists an integral surface $S$ through $x_{0}$.
\end{corollary}

\section{Global Frobenius Theorem\label{SectionDist&Foliations}}

The goal of this section is to recast Theorem \ref{FrobeniusThm} in the
language of distributions and foliations, and so we begin with several
definitions. $M$ is, as ever, a locally complete metric space.

\begin{definition}
A \textbf{distribution} $\Delta$ on $M$ is a set of curves in $M$.
\end{definition}

\begin{example}
A single arc field $X$ gives a distribution by forgetting $X$ is continuous in
its space variable $x$, and defining $\Delta=\left\{  X\left(  x,\cdot\right)
:x\in M\right\}  $. Any union of arc fields similarly gives a distribution.

Given two arc fields $X$ and $Y$, their \textbf{linear span} is a
distribution:%
\[
\Delta\left(  X,Y\right)  =\left\{  \left(  aX+bY\right)  \left(
x,\cdot\right)  :a,b\in\mathbb{R},x\in M\right\}  \text{.}%
\]
\end{example}

The direct sum of an arbitrary collection of arc fields similarly gives a
distribution, defined with finite summands. All of the following definitions
can, of course, be made with arbitrary indexing sets; but we will only use
finite sets of generators in the theorems of this paper.

Denote $\Delta_{x}:=\left\{  c\in\Delta:c\left(  0\right)  =x\right\}  $.

\begin{definition}
$X$ is $($locally uniformly$)$ \textbf{transverse} to $\Delta$ if for all
$x_{0}\in M$ there exists a $\delta>0$ such that for all $x\in B\left(
x_{0},\delta\right)  $ we have
\[
d\left(  X_{x}\left(  t\right)  ,c\left(  s\right)  \right)  \geq\delta\left(
\left|  s\right|  +\left|  t\right|  \right)
\]
for all $c\in\Delta_{x}$ and all $\left|  s\right|  ,\left|  t\right|
<\delta$. The arc fields $\overset{1}{X},$ $\overset{2}{X},$ $...,$
$\overset{n}{X}$ are \textbf{transverse} to each other if for each
$i\in\left\{  1,...,n\right\}  $ we have $\overset{i}{X}$ transverse to%
\[
\Delta\left(  \overset{1}{X},\overset{2}{X},...,\overset{i-1}{X},\overset
{i+1}{X},...,\overset{n}{X}\right)  \text{.}%
\]
\end{definition}

For $y\in M$ define%
\[
d\left(  y,\Delta_{x}\right)  :=\inf\left\{  d\left(  y,c\left(  t\right)
\right)  :c\in\Delta_{x}\text{ and }t\in dom\left(  c\right)  \right\}
\text{.}%
\]
If $\Delta_{x}=\emptyset$ then, as usual, the distance is $\infty$ by
definition. So if $X$ is transverse to $\Delta$ then if for all $x_{0}\in M$
there exists a $\delta>0$ such that for all $x\in B\left(  x_{0}%
,\delta\right)  $ we have
\[
d\left(  X_{x}\left(  t\right)  ,\Delta\right)  \geq\delta\left|  t\right|
\]
for all $\left|  t\right|  <\delta$.

\begin{definition}
$X$ is \textbf{tangent} to $\Delta$ if for each $x\in M$
\[
d\left(  X_{x}\left(  t\right)  ,\Delta_{x}\right)  =o\left(  t\right)
\text{.}%
\]
If this distance is $o\left(  t\right)  $ locally uniformly in $x\in M$ then
$X $ is \textbf{locally uniformly tangent} to $\Delta$, denoted $X\sim\Delta$.

Two distributions $\Delta$ and $\widetilde{\Delta}$ are \textbf{tangent} if
for each $c\in\Delta$ there exists $\widetilde{c}\in\widetilde{\Delta}$ such
that $\widetilde{c}$ is tangent to $c$ $($at $t=0),$ and \textit{vice-versa,}
for each $\widetilde{c}\in\widetilde{\Delta}$ there exists $c\in\Delta$ such
that $c$ is tangent to $\widetilde{c}$. Local uniform tangency is defined in
the obvious way, and denoted $\Delta\sim\widetilde{\Delta}$. Again, $\sim$ is
an equivalence relation.
\end{definition}

\begin{definition}
A distribution $\Delta$ is $n$-\textbf{dimensional} if there exists a set of
$n$ transverse arc fields $\left\{  \overset{1}{X},\overset{2}{X}%
,...,\overset{n}{X}\right\}  $ which all mutually close and satisfy E1 and E2
such that $\Delta\sim\Delta\left(  \overset{1}{X},\overset{2}{X}%
,...,\overset{n}{X}\right)  $.
\end{definition}

Given $X$, if there exist Lipschitz functions $a_{k}:M\rightarrow\mathbb{R}$
such that $X\sim\overset{n}{\underset{k=1}{\sum}}a_{k}\overset{k}{X}$ then
clearly $X\sim\Delta\left(  \overset{1}{X},\overset{2}{X},...,\overset{n}%
{X}\right)  $.

\begin{definition}
An $n$-dimensional distribution $\Delta\sim\Delta\left(  \overset{1}%
{X},\overset{2}{X},...,\overset{n}{X}\right)  $ is \textbf{involutive} if for
each choice of $i,j\in\left\{  1,...,n\right\}  $ we have%
\[
\left[  \overset{i}{X},\overset{j}{X}\right]  \sim\Delta\text{.}%
\]
\end{definition}

\begin{definition}
An \textbf{surface} $S$ in $M$ is an $n$-dimensional topological manifold
$S\subset M$. A surface is \textbf{locally uniformly tangent} to an arc field
$X$, denoted $X\sim S$, if $d\left(  X_{t}\left(  x\right)  ,S\right)
=o\left(  t\right)  $ locally uniformly in $x$.

A surface is said to be an \textbf{integral surface} for an $n$-dimensional
distribution $\Delta\sim\Delta\left(  \overset{1}{X},\overset{2}%
{X},...,\overset{n}{X}\right)  $ if $\overset{n}{\underset{k=1}{\sum}}%
a_{k}\overset{k}{X}\sim S$ for any choice of Lipschitz functions
$a_{k}:M\rightarrow\mathbb{R}$.

A distribution is said to be \textbf{integrable} if there exists an integral
surface through every point in $M$.
\end{definition}

Theorem \ref{FrobeniusThm} has the following corollary:

\begin{theorem}
\label{ThmInvol=>Integ}An $n$-dimensional involutive distribution is integrable.
\end{theorem}

\begin{proof}
$n=1$ is Theorem \ref{CL}. $n=2$ is Theorem \ref{FrobeniusThm}. Now proceed by
induction. We do enough of the case $n=3$ to suggest the path; and much of
this is copied from the proof of Theorem \ref{FrobeniusThm}.

Choose $x_{0}\in M$. Let $X,Y,$ and $Z$ be the transverse arc fields
guaranteed in the definition of a 3-dimensional distribution. If we find an
integral surface $S$ for $\Delta\left(  X,Y,Z\right)  $ through $x_{0}$ then
obviously $S$ is an integral surface for $\Delta$. Let $F,G$, and $H$ be the
local flows of $X,Y$, and $Z$ and define%
\[
S:=\left\{  F_{t}G_{s}H_{r}\left(  x_{0}\right)  :\left|  r\right|  ,\left|
s\right|  ,\left|  t\right|  <\delta\right\}
\]
with $\delta>0$ chosen small enough as in the proof of Theorem
\ref{FrobeniusThm} so that $S$ is a three dimensional manifold. Again we may
assume $\delta$ is also chosen small enough so that throughout $S$ the
functions $\left|  a_{k}\right|  $ are bounded by $A$, the constants $\Lambda
$, $\Omega$, and $\rho$ for $X,Y$ and $Z$ hold uniformly, and the closure of
$B\left(  x,3\delta\left(  \rho+1\right)  \right)  $ is complete. Notice
\[
\underline{S}:=\left\{  G_{s}H_{r}\left(  x_{0}\right)  :\left|  r\right|
,\left|  s\right|  <\delta\right\}
\]
is an integral surface through $x_{0}$ for $\Delta\left(  Y,Z\right)  $ by the
proof of Theorem \ref{FrobeniusThm}. Notice $S\sim X$ by construction, but it
is not immediately clear $S\sim a^{\prime}X+b^{\prime}Y+c^{\prime}Z $ for
arbitrarily chosen $a^{\prime},b^{\prime},c^{\prime}\in\mathbb{R}$. Again we
really only need to show $S\sim a^{\prime}F+b^{\prime}G+c^{\prime}H$ for an
arbitrary point $z:=F_{t}G_{s}H_{r}\left(  x_{0}\right)  \in S$, and again it
is sufficient to prove%
\[
\left(  F_{t}\right)  ^{\ast}\left(  a^{\prime}F+b^{\prime}G+c^{\prime
}H\right)  \sim S\text{\qquad at\qquad}y=G_{s}H_{r}\left(  x_{0}\right)
\]
by the construction of $S$. Continue as above adapting the same tricks from
the proof of Theorem \ref{FrobeniusThm} to the extra dimension.
\end{proof}

Similar to the definition for a surface, an arbitrary set $S\subset M$ is
defined to be \textbf{locally uniformly tangent} to $X$ if%
\[
d\left(  X_{t}\left(  y\right)  ,S\right)  =o\left(  t\right)
\]
locally uniformly for $y\in S$, denoted $S\sim X$.

\begin{lemma}
\label{LemmaNagumo}Let $\sigma_{x}:\left(  \alpha,\beta\right)  \rightarrow
U\subset M$ be a solution to $X$ which meets Condition E1 with uniform
constant $\Lambda$ on a neighborhood $U$. Assume $S\subset U$ is a closed set
with $S\sim X$. Then%
\[
d\left(  \sigma_{x}\left(  t\right)  ,S\right)  \leq e^{\Lambda\left|
t\right|  }d\left(  x,S\right)  \text{\ for all }t\in\left(  \alpha
,\beta\right)  \text{.}%
\]
\end{lemma}

\begin{proof}
(Adapted from the proof of Theorem \ref{ExpGrowth}.)

We check only $t>0$. Define%
\[
g\left(  t\right)  :=e^{-\Lambda t}d\left(  \sigma_{x}\left(  t\right)
,S\right)  \text{.}%
\]
For $h\geq0$, we have
\begin{align*}
&  g\left(  t+h\right)  -g\left(  t\right) \\
&  =e^{-\Lambda\left(  t+h\right)  }d\left(  \sigma_{x}\left(  t+h\right)
,S\right)  -e^{-\Lambda t}d\left(  \sigma_{x}\left(  t\right)  ,S\right) \\
&  \leq e^{-\Lambda\left(  t+h\right)  }\left[  d\left(  \sigma_{x}\left(
t+h\right)  ,X_{h}\left(  \sigma_{x}\left(  t\right)  \right)  \right)
+d\left(  X_{h}\left(  \sigma_{x}\left(  t\right)  \right)  ,X_{h}\left(
y\right)  \right)  +d\left(  X_{h}\left(  y\right)  ,S\right)  \right] \\
&  -e^{-\Lambda t}d\left(  \sigma_{x}\left(  t\right)  ,S\right)
\end{align*}
for any $y\in S,$ which in turn is%
\begin{align*}
&  \leq e^{-\Lambda\left(  t+h\right)  }\left[  d\left(  X_{h}\left(
\sigma_{x}\left(  t\right)  \right)  ,X_{h}\left(  y\right)  \right)
+o\left(  h\right)  \right]  -e^{-\Lambda t}d\left(  \sigma_{x}\left(
t\right)  ,S\right) \\
&  \leq e^{-\Lambda t}e^{-\Lambda h}d\left(  \sigma_{x}\left(  t\right)
,y\right)  \left(  1+\Lambda h\right)  -e^{-\Lambda t}d\left(  \sigma
_{x}\left(  t\right)  ,S\right)  +o\left(  h\right) \\
&  =\left[  e^{-\Lambda h}\left(  1+\Lambda h\right)  d\left(  \sigma
_{x}\left(  t\right)  ,y\right)  -d\left(  \sigma_{x}\left(  t\right)
,S\right)  \right]  e^{-\Lambda t}+o\left(  h\right)  \text{.}%
\end{align*}
Therefore%
\[
g\left(  t+h\right)  -g\left(  t\right)  \leq\left[  e^{-\Lambda h}\left(
1+\Lambda h\right)  -1\right]  e^{-\Lambda t}d\left(  \sigma_{x}\left(
t\right)  ,S\right)  +o\left(  h\right)
\]
since $y$ was arbitrary in $S$. Thus%
\begin{align*}
&  g\left(  t+h\right)  -g\left(  t\right) \\
&  \leq o\left(  h\right)  e^{-\Lambda t}d\left(  \sigma_{x}\left(  t\right)
,S\right)  +o\left(  h\right)  =o\left(  h\right)  \left(  g\left(  t\right)
+1\right)  .
\end{align*}
Hence, the upper forward derivative of $g\left(  t\right)  $ is nonpositive;
i.e.,
\[
D^{+}g\left(  t\right)  :=\overline{\lim_{h\rightarrow0^{+}}}\,\left(
\frac{g\left(  t+h\right)  -g\left(  t\right)  }{h}\right)  \leq0.
\]
Consequently, $g\left(  t\right)  \leq g\left(  0\right)  $ or
\[
d\left(  \sigma_{x}\left(  t\right)  ,S\right)  \leq e^{\Lambda t}d\left(
\sigma_{x}\left(  0\right)  ,S\right)  =e^{\Lambda t}d\left(  x,S\right)  .
\]
\end{proof}

Choosing $x\in S$ in Lemma \ref{LemmaNagumo} gives the following metric space
generalization of the Nagumo-Brezis Invariance Theorem (Example \ref{Banach
Example} shows how this generalizes the Banach space setting). We state and
prove only the bidirectional case; the case for forward flows is easily
adapted \textit{mutatis mutandis}. Cf. \cite{Motreanu Pavel} for an exposition
on general invariance theorems.

\begin{theorem}
\label{ThmNagumo}Let $X$ satisfy E1 and E2 and assume a closed set $S\subset
M$ has $S\sim X$. Then for any $x\in S$ we have $F_{t}\left(  x\right)  \in S$
for all $t\in\left(  \alpha_{x},\beta_{x}\right)  $. I.e., $S$ is an
\textbf{invariant set} under the flow $F$.
\end{theorem}

\begin{theorem}
\label{ThmUniqueIntSurfs}The integral surfaces guaranteed by Theorem
\ref{ThmInvol=>Integ} are unique in the following sense: if $S_{1}$ and
$S_{2}$ are integral surfaces through $x\in M$, then $S_{1}\cap S_{2}$ is an
integral surface.
\end{theorem}

\begin{proof}
The case $n=1$ is true by the uniqueness of integral curves.

For higher dimensions $n$, Theorem \ref{ThmNagumo} guarantees $S_{1}$ and
$S_{2}$ contain local integral curves for $\overset{n}{\underset{k=1}{\sum}%
}a_{k}\overset{k}{X}$ for all choices of $a_{k}\in\mathbb{R}$ with initial
condition $x$. Since the $\overset{k}{X}$ are transverse, there is a small
neighborhood of $x$ on which all the choices of the parameters $a_{k}$ give
local non-intersecting curves in $M$ which fill up $n$ dimensions.
\end{proof}

Therefore, by continuation we have a unique maximal connected integral surface
through each point.

\begin{definition}
A \textbf{foliation} partitions $M$ into a set of subsets $\Phi:=\left\{
\mathcal{L}_{i}\right\}  _{i\in I}$ for some indexing set $I$, where the
subsets $\mathcal{L}_{i}\subset M$ $\left(  \text{called \textbf{leaves}%
}\right)  $ are disjoint, connected topological manifolds each having the same dimension.

A foliation $\Phi$ is \textbf{tangent} to a distribution $\Delta$ if the
leaves are integral surfaces. When a foliation exists which is tangent to a
distribution $\Delta$ we say $\Delta$ \textbf{foliates} $M$.
\end{definition}

\begin{theorem}
\label{ThmInvol=>Foliat}An $n$-dimensional involutive distribution has a
unique tangent foliation.
\end{theorem}

\begin{proof}
Theorems \ref{ThmInvol=>Integ} and \ref{ThmUniqueIntSurfs} guarantee the
existence of the leaves, i.e., the unique maximal integral surfaces.
\end{proof}

The converse of this result is easy to prove in the classical context on a
Banach space. I do not believe it is true here. Instead we have the following
partial converse. Cf. Remark \ref{Rem2ndObracket}.

\begin{proposition}
Let $\Delta\sim\Delta\left(  \overset{1}{X},\overset{2}{X},...,\overset{n}%
{X}\right)  $ be an $n$-dimensional distribution with $\overset{i}{X}%
\approx\overset{i}{F}$ where $\overset{i}{F}$ is the local flow for
$\overset{i}{X}$. If $\Delta$ foliates $M$ then $\Delta$ is involutive.
\end{proposition}

\begin{proof}
Remark \ref{Rem2ndObracket} gives $\left[  \overset{i}{X},\overset{j}%
{X}\right]  \sim\left[  \overset{i}{F},\overset{j}{F}\right]  $ and Theorem
\ref{ThmNagumo} gives $\left[  \overset{i}{F},\overset{j}{F}\right]
_{t}\left(  x\right)  \in\mathcal{L}_{i}$ if $x\in\mathcal{L}_{i}$ so $\left[
\overset{i}{F},\overset{j}{F}\right]  \sim\Delta$.
\end{proof}

Collecting all these results we have the following version of the Global
Frobenius Theorem.

\begin{theorem}
\label{ThmGlobFrob}Let $\Delta\sim\Delta\left(  \overset{1}{X},\overset{2}%
{X},...,\overset{n}{X}\right)  $ be an $n$-dimensional distribution on a
locally complete metric space $M$, with $\overset{i}{X}\approx\overset{i}{F} $
where $\overset{i}{F}$ is the local flow for $\overset{i}{X}$. The following
are equivalent$:$

1. $\Delta$ is involutive

2. $\Delta$ is integrable

3. $\Delta$ foliates $M$.
\end{theorem}

\section{Commutativity of Flows\label{SectionCommut}}

\begin{theorem}
\label{ThmCommute}Assume $X$ and $Y$ satisfy E1 and E2 on a locally complete
metric space $M$. Let $F$ and $G$ be the local flows of $X$ and $Y$. Then
$\left[  F,G\right]  \sim0$ if and only if $F$ and $G$ commute, i.e.,%
\[
F_{t}G_{s}\left(  x\right)  =G_{s}F_{t}\left(  x\right)  \text{,\qquad
i.e.,\qquad}F_{t}^{\ast}\left(  G\right)  =G\text{.}%
\]
\end{theorem}

\begin{proof}
The assumption $\left[  F,G\right]  \sim aX+bY$ with $a=b=0$ allows us to copy
the approach in the proof of Theorem \ref{FrobeniusThm}. Let $\delta>0$ be
chosen small enough so

1. the functions $\left|  a\right|  $ and $\left|  b\right|  $ are bounded

2. the constants $\Lambda$, $\Omega$, and $\rho$ for $X$ and $Y$ hold uniformly

3. $\left[  F,G\right]  \sim0$ uniformly

\noindent all on $S:=B\left(  x,2\delta\left(  \rho+1\right)  \right)  $ and
that $S$ is also complete. We check $t>0$. Since $F_{t}^{\ast}\left(
G\right)  $ and $G$ are both local flows, we only need to show they are
tangent to each other and then they must be equal by uniqueness of solutions.

Imagine being in the context of differentiable manifolds. There, for vector
fields $f$ and $g$ with local flows $F$ and $G$, we would have%
\[
\underset{h\rightarrow0}{\lim}\frac{F_{h}^{\ast}\left(  g\right)  -g}%
{h}=\mathcal{L}_{f}g=\left[  f,g\right]  =0
\]
so $F_{h}^{\ast}\left(  g\right)  =g+o\left(  h\right)  $ and thus we expect%
\[
F_{h}^{\ast}\left(  g\right)  =g+o\left(  h\right)  \text{.}%
\]
We might use this idea as before with the linearity of pull-back (Proposition
\ref{PullbackLinear}) to get%
\[
F_{t}^{\ast}\left(  g\right)  =\underset{n\rightarrow\infty}{\lim}%
F_{t/n}^{\ast\left(  n\right)  }\left(  g\right)  =\underset{n\rightarrow
\infty}{\lim}g+no\left(  1/n\right)  =g
\]
as desired.

Now in our context of metric spaces with $t>0$, line $\left(
\text{\ref{LieD=BracketArcFieldversion}}\right)  $ again gives%
\[
F_{t/n}^{\ast}\left(  G\right)  _{t/n}\left(  x\right)  =\left(  \tfrac{t}%
{n}\left[  F,G\right]  +G\right)  _{t/n}\left(  x\right)  \text{.}%
\]
For $t<0$ one would use $\left(  \text{\ref{LieD=BracketArcFieldversion2}%
}\right)  $. Also we again have%
\begin{gather*}
\left[  F,G\right]  \sim0\qquad\text{implies}\\
d\left(  \left(  \tfrac{t}{n}\left[  F,G\right]  \right)  _{t/n}\left(
x\right)  ,x\right)  =o\left(  \tfrac{1}{n^{2}}\right)  \text{.}%
\end{gather*}
Using these tricks (and Theorem \ref{ExpGrowth} in the fourth line following)
gives%
\begin{align*}
&  d\left(  \left(  F_{t}^{\ast}\left(  G\right)  \right)  _{t/n}\left(
x\right)  ,G_{t/n}\left(  x\right)  \right)  =d\left(  \left(  F_{t/n}%
^{\ast\left(  n-1\right)  }F_{t/n}^{\ast}\left(  G\right)  \right)
_{t/n}\left(  x\right)  ,G_{t/n}\left(  x\right)  \right) \\
&  =d\left(  F_{t/n}^{\ast\left(  n-1\right)  }\left(  \tfrac{t}{n}\left[
F,G\right]  +G\right)  _{t/n}\left(  x\right)  ,G_{t/n}\left(  x\right)
\right) \\
&  \leq d\left(  F_{t/n}^{\ast\left(  n-1\right)  }\left(  G_{t/n}\tfrac{t}%
{n}\left[  F,G\right]  _{t/n}\left(  x\right)  \right)  ,F_{t/n}^{\ast\left(
n-1\right)  }G_{t/n}\left(  x\right)  \right)  +d\left(  F_{t/n}^{\ast\left(
n-1\right)  }G_{t/n}\left(  x\right)  ,G_{t/n}\left(  x\right)  \right) \\
&  \leq d\left(  G_{t/n}\tfrac{t}{n}\left[  F,G\right]  _{t/n}\left(
y\right)  ,G_{t/n}\left(  y\right)  \right)  e^{\Lambda_{X}\frac{t\left(
n-1\right)  }{n}}+d\left(  F_{t/n}^{\ast\left(  n-1\right)  }G_{t/n}\left(
x\right)  ,G_{t/n}\left(  x\right)  \right)
\end{align*}
where $y:=F_{\left(  n-1\right)  t/n}\left(  x\right)  $%
\[
\leq d\left(  \tfrac{t}{n}\left[  F,G\right]  _{t/n}\left(  y\right)
,y\right)  e^{\Lambda_{Y}t/n}e^{\Lambda_{X}\frac{t\left(  n-1\right)  }{n}%
}+d\left(  F_{t/n}^{\ast\left(  n-1\right)  }G_{t/n}\left(  x\right)
,G_{t/n}\left(  x\right)  \right)
\]
and so%
\begin{align*}
&  d\left(  \left(  F_{t}^{\ast}\left(  G\right)  \right)  _{t/n}\left(
x\right)  ,G_{t/n}\left(  x\right)  \right) \\
&  \leq d\left(  F_{t/n}^{\ast\left(  n-1\right)  }G_{t/n}\left(  x\right)
,G_{t/n}\left(  x\right)  \right)  +e^{\Lambda_{Y}t/n+\Lambda_{X}%
\frac{t\left(  n-1\right)  }{n}}o_{1}\left(  \tfrac{1}{n^{2}}\right)
\end{align*}
where $o_{1}\left(  \tfrac{1}{n^{2}}\right)  :=d\left(  \tfrac{t}{n}\left[
F,G\right]  _{t/n}\left(  y\right)  ,y\right)  $.

Iterating this result gives%
\begin{align*}
&  d\left(  \left(  F_{t/n}^{\ast n}\left(  G\right)  \right)  _{t/n}\left(
x\right)  ,G_{t/n}\left(  x\right)  \right) \\
&  \leq d\left(  F_{t/n}^{\ast\left(  n-1\right)  }G_{t/n}\left(  x\right)
,G_{t/n}\left(  x\right)  \right)  +e^{\Lambda_{Y}t/n+\Lambda_{X}%
\frac{t\left(  n-1\right)  }{n}}o_{1}\left(  \tfrac{1}{n^{2}}\right) \\
&  \leq d\left(  F_{t/n}^{\ast\left(  n-2\right)  }G_{t/n}\left(  x\right)
,G_{t/n}\left(  x\right)  \right)  +e^{\Lambda_{Y}t/n+\Lambda_{X}%
\frac{t\left(  n-2\right)  }{n}}o_{2}\left(  \tfrac{1}{n^{2}}\right)
+e^{\Lambda_{Y}t/n+\Lambda_{X}\frac{t\left(  n-1\right)  }{n}}o_{1}\left(
\tfrac{1}{n^{2}}\right) \\
&  \leq...\leq d\left(  F_{t/n}^{0}G_{t/n}\left(  x\right)  ,G_{t/n}\left(
x\right)  \right)  +e^{\Lambda_{Y}t/n}\overset{n}{\underset{i=1}{%
%TCIMACRO{\tsum }%
%BeginExpansion
{\textstyle\sum}
%EndExpansion
}}o_{i}\left(  \tfrac{1}{n^{2}}\right)  e^{\Lambda_{X}\frac{t\left(
n-i\right)  }{n}}\\
&  =e^{\Lambda_{Y}t/n}\overset{n}{\underset{i=1}{%
%TCIMACRO{\tsum }%
%BeginExpansion
{\textstyle\sum}
%EndExpansion
}}o_{i}\left(  \tfrac{1}{n^{2}}\right)  e^{\Lambda_{X}\frac{t\left(
n-i\right)  }{n}}%
\end{align*}
where $o_{i}\left(  \tfrac{1}{n^{2}}\right)  :=d\left(  \tfrac{t}{n}\left[
F,G\right]  _{t/n}\left(  y_{i}\right)  ,y_{i}\right)  $ and $y_{i}%
:=F_{\left(  n-i\right)  t/n}\left(  x\right)  $. Since $d\left(  \tfrac{t}%
{n}\left[  F,G\right]  _{t/n}\left(  y\right)  ,y\right)  =o\left(  \tfrac
{1}{n^{2}}\right)  $ uniformly for $y\in B\left(  x,2\delta\left(
\rho+1\right)  \right)  $ we have%
\begin{align*}
&  d\left(  \left(  F_{t/n}^{\ast n}\left(  G\right)  \right)  _{t/n}\left(
x\right)  ,G_{t/n}\left(  x\right)  \right) \\
&  \leq e^{\Lambda_{Y}t/n}\overset{n}{\underset{i=1}{%
%TCIMACRO{\tsum }%
%BeginExpansion
{\textstyle\sum}
%EndExpansion
}}o_{i}\left(  \tfrac{1}{n^{2}}\right)  e^{\Lambda_{X}\frac{t\left(
n-i\right)  }{n}}=o\left(  \tfrac{1}{n^{2}}\right)  e^{\Lambda_{Y}t/n}%
\overset{n}{\underset{i=1}{%
%TCIMACRO{\tsum }%
%BeginExpansion
{\textstyle\sum}
%EndExpansion
}}e^{\Lambda_{X}\frac{t\left(  n-i\right)  }{n}}\\
&  =o\left(  \tfrac{1}{n^{2}}\right)  e^{\Lambda_{Y}t/n}e^{\Lambda_{X}%
t}\overset{n}{\underset{i=1}{%
%TCIMACRO{\tsum }%
%BeginExpansion
{\textstyle\sum}
%EndExpansion
}}\left(  e^{-\frac{t}{n}}\right)  ^{i}=o\left(  \tfrac{1}{n^{2}}\right)
e^{\Lambda_{Y}t/n+\Lambda_{X}t}\frac{1-\left(  e^{-\frac{t}{n}}\right)
^{n+1}}{1-\left(  e^{-\frac{t}{n}}\right)  }\text{.}%
\end{align*}
So%
\[
d\left(  \left(  F_{t}^{\ast}\left(  G\right)  \right)  _{t/n}\left(
x\right)  ,G_{t/n}\left(  x\right)  \right)  =o\left(  \tfrac{1}{n}\right)
\]
and $F_{t}^{\ast}\left(  G\right)  \sim G$ by the same argument at the last
paragraph of the proof of Theorem \ref{FrobeniusThm}.

The converse is trivial.
\end{proof}

Using Example \ref{Banach Example}, this theorem applies to the non-locally
compact setting with nonsmooth vector fields. \cite{RampSuss}, another paper
which inspires this monograph, obtains similar results with a very different approach.

\section{Examples\label{SectionExs}}

\begin{example}
Let $M$ be a Banach space. First let $X$ and $Y$ be translations in the
directions of $u$ and $v\in M$%
\[
X_{t}\left(  x\right)  :=x+tu\text{\qquad}Y_{t}\left(  x\right)  :=x+tv
\]
then $F=X$ and $G=Y$ for $\left|  t\right|  \leq1$. Obviously $\left[
F,G\right]  =0$ and the flows commute.

Next consider the dilations $X$ and $Y$ about $u$ and $v\in M$%
\[
X_{t}\left(  x\right)  :=\left(  1+t\right)  \left(  x-u\right)
+u\text{\qquad}Y_{t}\left(  x\right)  :=\left(  1+t\right)  \left(
x-v\right)  +v\text{.}%
\]
The flows are computable with little effort using Euler curves, e.g.,%
\[
F_{t}\left(  x\right)  =\underset{n\rightarrow\infty}{\lim}X_{t/n}^{\left(
n\right)  }\left(  x\right)  =e^{t}x-\left(  e^{t}-1\right)  u\text{.}%
\]
Then for $t\geq0$%
\begin{align*}
&  \left[  F,G\right]  _{t^{2}}\left(  x\right) \\
&  =G_{-t}F_{-t}G_{t}F_{t}\left(  x\right) \\
&  =e^{-t}\left[  e^{-t}\left(  e^{t}\left[  e^{t}x-\left(  e^{t}-1\right)
u\right]  -\left(  e^{t}-1\right)  v\right)  -\left(  e^{-t}-1\right)
u\right]  -\left(  e^{-t}-1\right)  v\\
&  =x-u+e^{-t}u-e^{-t}v+e^{-2t}v-e^{-2t}u+e^{-t}u-e^{-t}v+v\\
&  =x+\left(  v-u\right)  \left(  e^{-t}-1\right)  ^{2}%
\end{align*}
so $\left[  F,G\right]  \sim Z$ where $Z$ is the translation $Z_{t}\left(
x\right)  :=x+t\left(  v-u\right)  $ since, for instance with $t>0$%
\begin{align*}
&  d\left(  \left[  F,G\right]  _{t}\left(  x\right)  ,Z_{t}\left(  x\right)
\right) \\
&  =\left|  v-u\right|  \left|  \left(  e^{-\sqrt{t}}-1\right)  ^{2}-t\right|
=\left|  t\right|  \left|  v-u\right|  \left|  \left(  \frac{e^{-\sqrt{t}}%
-1}{\sqrt{t}}\right)  ^{2}-1\right|  =o\left(  t\right)  \text{.}%
\end{align*}
Hence the distribution $\Delta\left(  X,Y\right)  $ is not involutive.
However, this shows three dilations $X,Y,Z$ about linearly independent $u,v,w
$ generate all translations using their brackets. Using the same tricks we've
just employed, it is easy to check the bracket of a dilation and a translation
is tangent to a translation, e.g., if $F_{t}\left(  x\right)  :=x+tu$ and
$G_{t}\left(  x\right)  :=e^{t}x$ $($dilation about $0)$ then $\left[
F,G\right]  \sim F$ since for $t>0$%
\[
\left[  F,G\right]  _{t^{2}}\left(  x\right)  =G_{-t}F_{-t}G_{t}F_{t}\left(
x\right)  =e^{-t}\left[  e^{t}\left[  x+tu\right]  -tu\right]  =x+tu\left(
1-e^{-t}\right)
\]
and so
\[
d\left(  \left[  F,G\right]  _{t}\left(  x\right)  ,F_{t}\left(  x\right)
\right)  =\left|  tu\right|  \left|  \tfrac{1-e^{-\sqrt{t}}}{\sqrt{t}%
}-1\right|  =o\left(  t\right)  \text{.}%
\]
To summarize:%
\begin{equation}%
\begin{array}
[c]{ll}%
\Delta\left(  translations\right)  & involutive\\
\Delta\left(  dilations\right)  & \text{not involutive}\\
\Delta\left(  dilations,translations\right)  & involutive\text{.}%
\end{array}
\label{ExSummary}%
\end{equation}
\end{example}

The previous example holds with minor modifications on the metric space
$\left(  H\left(  \mathbb{R}^{n}\right)  ,d_{H}\right)  $ where $H\left(
\mathbb{R}^{n}\right)  $ is the set of non-void compact subsets of
$\mathbb{R}^{n}$ and $d_{H}$ is the Hausdorff metric. Theorem
\ref{ThmInvol=>Foliat} gives foliations.

\begin{example}
[two parameter decomposition of $L^{2}$]\label{ExL2decomp}Now let $M$ be real
Hilbert space $L^{2}\left(  \mathbb{R}\right)  $. Since $M$ is Banach the
results of the previous example hold. Denote translation by the function $h\in
L^{2}\left(  \mathbb{R}\right)  $ by%
\[
X_{t}\left(  f\right)  :=f+th\text{.}%
\]
Now however, there is another obvious candidate for an elementary flow:
translation with respect to the variable $x$, i.e.,%
\[
Y_{t}\left(  f\right)  \left(  x\right)  :=f\left(  x+t\right)  \text{.}%
\]
Unlike dilation and translation, the dynamic engendered by $Y$ seemingly has
nothing to do with the vector space structure of $L^{2}\left(  \mathbb{R}%
\right)  $. In fact, despite appearances, $Y$ is a nonsmooth flow: notice for
example with the characteristic function $\chi$ as initial condition,%
\[
\left.  \frac{d}{dt}Y_{t}\left(  \chi_{\left[  0,1\right]  }\right)  \right|
_{t=0}\notin L^{2}\left(  \mathbb{R}\right)  \text{.}%
\]
Interpreted as a flow on a metric space, however, this is no obstacle. We
refer to $X$ as \textbf{vector space translation} and $Y$ as \textbf{function
translation}. Notice $X$ and $Y$ are their own flows $($for $\left|  t\right|
\leq1)$. It is straightforward to check $X$ \& $Y$ close when, for example,
$h\in C^{1}\left(  \mathbb{R}\right)  $ with derivative $h^{\prime}\in
L^{2}\left(  \mathbb{R}\right)  :$%
\begin{align*}
&  d\left(  Y_{s}X_{t}\left(  f\right)  ,X_{t}Y_{s}\left(  f\right)  \right)
\\
&  =\sqrt{\int\left(  f\left(  x+s\right)  +th\left(  x+s\right)  -\left[
f\left(  x+s\right)  +th\left(  x\right)  \right]  \right)  ^{2}dx}\\
&  =\left|  st\right|  \sqrt{\int\left(  \frac{h\left(  x+s\right)  -h\left(
x\right)  }{s}\right)  ^{2}dx}\\
&  =O\left(  \left|  st\right|  \right)
\end{align*}
uniformly. Since they obviously satisfy E1 and E2, Theorem \ref{ThmaX+bYflow}
promises a unique flow for their sum. This was introduced by Colombo and Corli
in \cite[section 5.2]{Colombo} with other interesting function space examples,
which they also characterize with partial differential equations.

Let us now compute the bracket. We check $t>0$ explicitly, skipping the case
$t\leq0$ though this is just as easy.%
\begin{align*}
&  \left[  X,Y\right]  _{t^{2}}\left(  f\right)  \left(  x\right) \\
&  =Y_{-t}X_{-t}Y_{t}X_{t}\left(  f\right)  \left(  x\right)  =Y_{-t}%
X_{-t}\left[  f\left(  x+t\right)  +th\left(  x+t\right)  \right] \\
&  =f\left(  x\right)  +th\left(  x\right)  -th\left(  x-t\right)  =f\left(
x\right)  +t^{2}\left[  \frac{h\left(  x\right)  -h\left(  x-t\right)  }%
{t}\right]  \text{.}%
\end{align*}
Defining a new arc field $Z_{t}\left(  f\right)  :=f+th^{\prime}$ we therefore
have%
\[
d\left(  \left[  X,Y\right]  _{t}\left(  f\right)  ,Z_{t}\left(  f\right)
\right)  =\left|  t\right|  \sqrt{\underset{\mathbb{R}}{\int}\left(
\frac{h\left(  x\right)  -h\left(  x-t\right)  }{t}-h^{\prime}\left(
x\right)  \right)  ^{2}dx}=o\left(  t\right)
\]
when $h\in C^{1}\left(  \mathbb{R}\right)  $ with $h^{\prime}\in L^{2}\left(
\mathbb{R}\right)  $. Thus $\left[  X,Y\right]  \sim Z$.

This has remarkable consequences. Using the idea of Chow's Theorem from
control theory $\left(  \text{also called the Chow-Rashevsky Theorem or
Hermes' Theorem}\right)  $, if the $\left(  n+1\right)  $-st derivative
$h^{\left[  n+1\right]  }$ is not contained in $span\left\{  h^{\left[
i\right]  }:0\leq i\leq n\right\}  $ then iterating the process of bracketing
$X$ and $Y$ generates a large space reachable via repeated compositions of $X$
and $Y$. Denoting%
\begin{equation}
\overset{n}{Z}_{t}\left(  f\right)  :=f+th^{\left[  n\right]  }
\label{ExL2dec10}%
\end{equation}
successive brackets of $X$ and $Y$ are%
\begin{gather}
\left[  X,Y\right]  \sim Z=:\overset{1}{Z}\nonumber\\
\left[  X\overset{2}{,}Y\right]  :=\left[  \left[  X,Y\right]  ,Y\right]
\sim\overset{2}{Z}\nonumber\\
\left[  X\overset{n}{,}Y\right]  :=\underset{n\text{ times}}{\underbrace
{\left[  \left[  ...\left[  \left[  X,Y\right]  ,Y\right]  ,...,Y\right]
,Y\right]  }}\sim\overset{n}{Z}\text{.} \label{ExL2dec5}%
\end{gather}
For notational purposes we set $\left[  X\overset{0}{,}Y\right]  :=X$. In the
particular case $h\left(  x\right)  :=e^{-x^{2}}$ all of $L^{2}\left(
\mathbb{R}\right)  $ is reachable by $X$ and $Y$.

To see this we apply the theory of orthogonal functions with the
Hermite\footnote{We may of course use other orthogonal families with a
different choice of $h$, particularly when the domain of interest is other
than $\mathbb{R}$; e.g., scaled Chebyshev polynomials for $\left[
0,2\pi\right)  $, etc. We expect many choices of $h$ give controllable systems
whether the brackets generate orthogonal sets or not.} polynomials%
\[
H_{n}\left(  x\right)  :=\left(  -1\right)  ^{n}e^{x^{2}}\frac{d^{n}}{dx^{n}%
}e^{-x^{2}}=\left(  -1\right)  ^{n}e^{x^{2}}h^{\left[  n\right]  }\left(
x\right)
\]
which have dense span in $L^{2}\left(  \mathbb{R}\right)  $ when multiplied by
$e^{-x^{2}/2}$. Those familiar with orthogonal expansions can predict the
rest; we review some of the details.%
\[
\left\{  \dfrac{1}{\sqrt{n!2^{n}\sqrt{\pi}}}H_{n}\left(  x\right)
e^{-x^{2}/2}:n\in\mathbb{N}\right\}
\]
is a basis of $L^{2}\left(  \mathbb{R}\right)  $ and is orthonormal since%
\begin{equation}
\underset{\mathbb{R}}{%
%TCIMACRO{\tint }%
%BeginExpansion
{\textstyle\int}
%EndExpansion
}H_{m}\left(  x\right)  H_{n}\left(  x\right)  e^{-x^{2}}dx=n!2^{n}\sqrt{\pi
}\delta_{mn}\text{.} \label{ExHermite2}%
\end{equation}
The Hermite polynomials also satisfy some useful relations%
\begin{equation}
H_{n+1}\left(  x\right)  =2xH_{n}\left(  x\right)  -2nH_{n-1}\left(  x\right)
\qquad\text{and}\qquad H_{n}^{\prime}\left(  x\right)  =2nH_{n-1}\left(
x\right)  \text{.} \label{ExHermite6}%
\end{equation}
Given any $g\in L^{2}\left(  \mathbb{R}\right)  $ it is possible to write%
\begin{equation}
g\left(  x\right)  =\overset{\infty}{\underset{n=0}{%
%TCIMACRO{\tsum }%
%BeginExpansion
{\textstyle\sum}
%EndExpansion
}}a_{n}\tfrac{1}{\sqrt{n!2^{n}\sqrt{\pi}}}H_{n}\left(  x\right)  e^{-x^{2}/2}
\label{ExHermite10}%
\end{equation}
$($equality in the $L^{2}$ sense$)$ where%
\[
a_{n}:=\tfrac{1}{\sqrt{n!2^{n}\sqrt{\pi}}}\underset{\mathbb{R}}{%
%TCIMACRO{\tint }%
%BeginExpansion
{\textstyle\int}
%EndExpansion
}g\left(  x\right)  H_{n}\left(  x\right)  e^{-x^{2}/2}\left(  x\right)
dx\in\mathbb{R}\text{.}%
\]
The necessity of this formula for $a_{n}$ can easily be checked by multiplying
both sides of $\left(  \ref{ExHermite10}\right)  $ by $H_{n}\left(  x\right)
e^{-x^{2}/2}$, integrating and applying $\left(  \ref{ExHermite2}\right)  $.
However, we want%
\[
g=\overset{\infty}{\underset{n=0}{%
%TCIMACRO{\tsum }%
%BeginExpansion
{\textstyle\sum}
%EndExpansion
}}c_{n}h^{\left[  n\right]  }%
\]
so apply the above process to $g\left(  x\right)  e^{x^{2}/2}$
instead\footnote{The function $g\left(  x\right)  e^{x^{2}/2}$ is no longer
necessarily $L^{2}$, of course, but here we lapse into the habit of ignoring
convergence issues as they are important for the theoretical proof that all of
$L^{2}\left(  \mathbb{R}\right)  $ is reachable with $X$ and $Y$, but not
central to this demonstration. This theoretical lapse is easily remedied by
multiplying by the characteristic function $\chi_{\left[  -m,m\right]  }$ to
guarantee all of the following integrals converge, then letting $m\rightarrow
\infty$ at the end.}. Then%
\begin{align*}
g\left(  x\right)  e^{x^{2}/2}  &  =\overset{\infty}{\underset{n=0}{%
%TCIMACRO{\tsum }%
%BeginExpansion
{\textstyle\sum}
%EndExpansion
}}b_{n}\tfrac{1}{\sqrt{n!2^{n}\sqrt{\pi}}}H_{n}\left(  x\right)  e^{-x^{2}%
/2}\text{\qquad so}\\
g  &  =\overset{\infty}{\underset{n=0}{%
%TCIMACRO{\tsum }%
%BeginExpansion
{\textstyle\sum}
%EndExpansion
}}c_{n}h^{\left[  n\right]  }%
\end{align*}
where%
\begin{align*}
b_{n}  &  :=\tfrac{1}{\sqrt{n!2^{n}\sqrt{\pi}}}\underset{\mathbb{R}}{%
%TCIMACRO{\tint }%
%BeginExpansion
{\textstyle\int}
%EndExpansion
}g\left(  x\right)  e^{x^{2}/2}H_{n}\left(  x\right)  e^{-x^{2}/2}\left(
x\right)  dx\text{\qquad so that}\\
c_{n}  &  :=\tfrac{\left(  -1\right)  ^{n}}{n!2^{n}\sqrt{2\pi}}\underset
{\mathbb{R}}{%
%TCIMACRO{\tint }%
%BeginExpansion
{\textstyle\int}
%EndExpansion
}g\left(  x\right)  h^{\left[  n\right]  }\left(  x\right)  e^{x^{2}}\left(
x\right)  dx\text{.}%
\end{align*}
Therefore when $N$ is large, $g$ is approximated by
\[
\overset{N}{\underset{n=0}{%
%TCIMACRO{\tsum }%
%BeginExpansion
{\textstyle\sum}
%EndExpansion
}}c_{n}h^{\left[  n\right]  }=F_{1}\left(  0\right)
\]
where $F$ is the flow of the arc field%
\[
\widetilde{X}:=\overset{N}{\underset{n=0}{%
%TCIMACRO{\tsum }%
%BeginExpansion
{\textstyle\sum}
%EndExpansion
}}c_{n}\left[  X\overset{n}{,}Y\right]
\]
which we follow for unit time starting with initial condition $0\in
L^{2}\left(  \mathbb{R}\right)  $. $F$ can, of course, be approximated by
Euler curves%
\[
F_{1}\left(  0\right)  =\underset{n\rightarrow\infty}{\lim}\widetilde{X}%
_{1/n}^{\left(  n\right)  }\left(  0\right)
\]
and since $\widetilde{X}$ is merely a $($complicated$)$ composition of $X$ and
$Y$, this gives us a simple algorithm for approximating any function $g$ with
only two simple flows.

Let us compute a basic example to illustrate this surprising fact. Choosing at
random $g\left(  x\right)  :=\chi_{\left[  0,1\right]  }\left(  x\right)  $,
the characteristic function of the unit interval, we have%
\begin{align*}
c_{n}  &  :=\tfrac{\left(  -1\right)  ^{n}}{n!2^{n}\sqrt{2\pi}}\underset
{0}{\overset{1}{%
%TCIMACRO{\tint }%
%BeginExpansion
{\textstyle\int}
%EndExpansion
}}H_{n}\left(  x\right)  dx=\tfrac{\left(  -1\right)  ^{n}}{2\left(
n+1\right)  n!2^{n}\sqrt{2\pi}}\left[  H_{n+1}\left(  1\right)  -H_{n+1}%
\left(  0\right)  \right]  \qquad\text{so, e.g.,}\\
c_{0}  &  =\tfrac{1}{\sqrt{2\pi}},\text{ }c_{1}=\tfrac{-1}{2\sqrt{2\pi}%
},\text{ }c_{2}=\tfrac{1}{12\sqrt{2\pi}},\text{ }c_{3}=\tfrac{1}{12\sqrt{2\pi
}},\text{ }c_{4}=\tfrac{1}{480\sqrt{2\pi}},\text{ etc.}%
\end{align*}
by $\left(  \ref{ExHermite6}\right)  $. Then stopping for the purposes of
illustration at $N=3$ our function $g$ is approximated by%
\[
\overset{3}{\underset{n=0}{%
%TCIMACRO{\tsum }%
%BeginExpansion
{\textstyle\sum}
%EndExpansion
}}c_{n}h^{\left[  n\right]  }\text{.}%
\]
Notice the flow of $\overset{i}{Z}$ from $\left(  \text{\ref{ExL2dec10}%
}\right)  $ is locally the same as $\overset{i}{Z}$ since it is just vector
space translation, so we will use the same symbol. All vector space
translations commute under $($arc field$)$ addition, and the arc field%
\[
\widetilde{Z}_{t}\left(  f\right)  :=\left(  c_{0}\overset{0}{Z}+c_{1}%
\overset{1}{Z}+c_{2}\overset{2}{Z}+c_{3}\overset{3}{Z}\right)  _{t}\left(
f\right)
\]
is locally equal to its flow. Obviously%
\[
\widetilde{Z}_{1}\left(  0\right)  =\overset{3}{\underset{n=0}{%
%TCIMACRO{\tsum }%
%BeginExpansion
{\textstyle\sum}
%EndExpansion
}}c_{n}h^{\left[  n\right]  }%
\]
and $\widetilde{Z}\sim\widetilde{X}$ where%
\begin{align*}
\widetilde{X}_{t}\left(  f\right)   &  :=\left(  c_{0}X+c_{1}\left[
X,Y\right]  +c_{2}\left[  \left[  X,Y\right]  ,Y\right]  +c_{3}\left[  \left[
\left[  X,Y\right]  ,Y\right]  ,Y\right]  \right)  _{t}\left(  f\right) \\
&  =\left(  \overset{3}{\underset{n=0}{%
%TCIMACRO{\tsum }%
%BeginExpansion
{\textstyle\sum}
%EndExpansion
}}c_{n}\left[  X\overset{n}{,}Y\right]  \right)  _{t}\left(  f\right)
\text{.}%
\end{align*}
Remember the arc field bracket and the arc field sum are defined as nothing
more than compositions of arc fields, e.g.,%
\begin{align*}
&  \left(  c_{0}X+c_{1}\left[  X,Y\right]  +c_{2}\left[  \left[  X,Y\right]
,Y\right]  \right)  _{t}\\
&  =\left[  \left[  X,Y\right]  ,Y\right]  _{c_{2}t}\left[  X,Y\right]
_{c_{1}t}X_{c_{0}t}%
\end{align*}
and, e.g., when $t>0$%
\begin{align*}
&  c_{2}\left[  \left[  X,Y\right]  ,Y\right]  _{t}\\
&  =Y_{-\sqrt{c_{2}t}}\left(  X_{-\sqrt{c_{2}t}}Y_{-\sqrt{c_{2}t}}%
X_{\sqrt{c_{2}t}}Y_{\sqrt{c_{2}t}}\right)  Y_{\sqrt{c_{2}t}}(Y_{-\sqrt{c_{2}%
t}}X_{-\sqrt{c_{2}t}}Y_{\sqrt{c_{2}t}}X_{\sqrt{c_{2}t}})\text{.}%
\end{align*}
Therefore this approximation of $g$ is achieved by computing the Euler curves
for $\widetilde{X}$ which is a complicated process $($with a simple formula$)$
of composing the elementary operations of function translation $(Y)$ and
vector space translation by the Gaussian $(X)$.

Continuing the example, for choices of $h$ other than the Gaussian it may be
the case that $h^{\left[  n+1\right]  }\in span\left\{  h^{\left[  i\right]
}:0\leq i\leq n\right\}  $. Then the space reachable by $X$ and $Y$ is
precisely limited. E.g., when $h$ is a trigonometric function from the
orthogonal Fourier decomposition of $L^{2}$ the parameter space is
two-dimensional, or when $h$ is an $n$-th order polynomial in the context of
$M=L^{2}\left[  a,b\right]  $ then the parameter space is $\left(  n+1\right)
$-dimensional.

Restating these results in different terminology: Controlling amplitude and
phase the 2-parameter system is holonomically constrained. Controlling phase
and superposition perturbation $(Y$ and $X)$ generates a larger space of
signals; how much $Y$ and $X$ deviate from holonomy depends on the choice of
perturbation function $h$. Consequently, a result for signal analysis is:
controlling two parameters is enough to generate any signal.
\end{example}

We collect some of the results of the previous example. Denote the
\textbf{reachable set} of $X$ and $Y$ by%
\[
R\left(  X,Y\right)  :=\left\{  Y_{s_{n}}X_{t_{n}}Y_{s_{n-1}}X_{t_{n-1}%
}...Y_{s_{1}}X_{t_{1}}\left(  0\right)  \in L^{2}\left(  \mathbb{R}\right)
:s_{i},t_{i}\in\mathbb{R},\text{ }n\in\mathbb{N}\right\}
\]
where $0\in L^{2}\left(  \mathbb{R}\right)  $ is the constant function.
$R\left(  X,Y\right)  $ is the set of all finite compositions of $X$ and $Y$.

\begin{theorem}
Let $h\in L^{2}\left(  \mathbb{R}\right)  $ be the Gaussian $h\left(
x\right)  :=e^{-x^{2}}$ and define
\[
X_{t}\left(  f\right)  :=f+th\qquad\text{and}\qquad Y_{t}\left(  f\right)
\left(  x\right)  :=f\left(  x+t\right)  \text{.}%
\]
Then $R\left(  X,Y\right)  $ is dense in $L^{2}\left(  \mathbb{R}\right)  $.
\end{theorem}

\begin{algorithm}
Let $g\in L^{2}\left(  \mathbb{R}\right)  $ be such that $\underset
{\mathbb{R}}{\int}\left[  g\left(  x\right)  e^{x^{2}/2}\right]  ^{2}%
dx<\infty$. Then
\[
g=\underset{n\rightarrow\infty}{\lim}\widetilde{X}_{1/n}^{\left(  n\right)
}\left(  0\right)
\]
where%
\begin{gather*}
\widetilde{X}:=\overset{\infty}{\underset{n=0}{%
%TCIMACRO{\tsum }%
%BeginExpansion
{\textstyle\sum}
%EndExpansion
}}c_{n}\left[  X\overset{n}{,}Y\right]  \text{\qquad with}\qquad c_{n}%
:=\tfrac{\left(  -1\right)  ^{n}}{n!2^{n}\sqrt{2\pi}}\underset{\mathbb{R}}{%
%TCIMACRO{\tint }%
%BeginExpansion
{\textstyle\int}
%EndExpansion
}g\left(  x\right)  h^{\left[  n\right]  }\left(  x\right)  e^{x^{2}}\left(
x\right)  dx\\
\text{and}\qquad\left[  X\overset{n}{,}Y\right]  :=\underset{n\text{ times}%
}{\underbrace{\left[  \left[  ...\left[  \left[  X,Y\right]  ,Y\right]
,...,Y\right]  ,Y\right]  }}%
\end{gather*}
and%
\[
\left[  X,Y\right]  \left(  f,t\right)  :=\left\{
\begin{array}
[c]{c}%
Y_{-\sqrt{t}}X_{-\sqrt{t}}Y_{\sqrt{t}}X_{\sqrt{t}}\left(  f\right) \\
X_{-\sqrt{\left|  t\right|  }}Y_{-\sqrt{\left|  t\right|  }}X_{\sqrt{\left|
t\right|  }}Y_{\sqrt{\left|  t\right|  }}\left(  f\right)
\end{array}
\right.
\begin{array}
[c]{c}%
\text{for }t\geq0\\
\text{for }t<0
\end{array}
\]
for any $f\in L^{2}\left(  \mathbb{R}\right)  $.
\end{algorithm}

\begin{example}
Let us continue Example \ref{ExL2decomp} with $M=L^{2}\left(  \mathbb{R}%
\right)  $ and%
\[
X_{t}\left(  f\right)  :=f+th\qquad\text{and}\qquad Y_{t}\left(  f\right)
\left(  x\right)  :=f\left(  x+t\right)
\]
which are vector space translation and function translation. Define the arc
fields%
\[
V_{t}\left(  f\right)  :=e^{t}f\qquad\text{and}\qquad W_{t}\left(  f\right)
\left(  x\right)  :=f\left(  e^{t}x\right)
\]
which may be thought of as \textbf{vector space dilation} $($about the point
$0\in M)$ and \textbf{function dilation} $($about the point $0\in\mathbb{R})
$. Again, $V$ and $W$ are coincident with their own flows. Using the same
approach as in Example \ref{ExL2decomp} it is easy to check the brackets
satisfy%
\[%
\begin{array}
[c]{ll}%
\left[  X,Y\right]  _{t}\left(  f\right)  =f+th^{\prime}+o\left(  t\right)  &
\qquad\left[  X,V\right]  _{t}\left(  f\right)  =f+th+o\left(  t\right) \\
\left[  X,W\right]  _{t}\left(  f\right)  \left(  x\right)  =f\left(
x\right)  +txh^{\prime}\left(  x\right)  +o\left(  t\right)  & \qquad\left[
Y,V\right]  =0\\
\left[  Y,W\right]  _{t}\left(  f\right)  \left(  x\right)  =f\left(
x-t\right)  +o\left(  t\right)  & \qquad\left[  V,W\right]  =0
\end{array}
\]
assuming for the $\left[  X,Y\right]  $ and $\left[  X,W\right]  $
calculations that $h\in C^{1}\left(  \mathbb{R}\right)  $ and $h^{\prime}\in
L^{2}\left(  \mathbb{R}\right)  $. Consequently%
\[%
\begin{array}
[c]{ll}%
\Delta\left(  X,Y\right)  & \text{may be highly non-involutive depending on
}h\text{,}\\
\Delta\left(  X,V\right)  & \text{is involutive, but }X\text{ and }V\text{ do
\textbf{not} commute,}\\
\Delta\left(  X,W\right)  & \text{may be highly non-involutive depending on
}h\text{,}\\
\Delta\left(  Y,V\right)  & \text{is involutive; }Y\text{ and }V\text{
commute,}\\
\Delta\left(  Y,W\right)  & \text{is involutive, but }Y\text{ and }W\text{ do
\textbf{not} commute,}\\
\Delta\left(  V,W\right)  & \text{is involutive; }V\text{ and }W\text{
commute.}%
\end{array}
\]
When $h$ is chosen correctly, $X$ and $W$ control many function spaces,
similarly to $X$ and $Y$. The four involutive distributions foliate
$L^{2}\left(  \mathbb{R}\right)  $.
\end{example}

\end{document}